\documentclass[letterpaper]{article}
\usepackage{aaai}
\usepackage{graphicx}
\usepackage{amssymb}
\usepackage{amsmath}
\usepackage{times}
\usepackage{enumerate}
\usepackage{multirow}
\usepackage{subfigure}
\usepackage{bm}
\usepackage[framed,numbered]{mcode}
\usepackage{mathtools}
\usepackage{amsthm}
\usepackage{algorithm2e}
\include{standardHeader}
\renewcommand{\a}{\mathbf{a}}
\newcommand{\y}{\mathbf{y}}

\newcommand{\ykk}{\mathbf{y}^{k+1}}
\newcommand{\x}{\mathbf{x}}
\newcommand{\xk}{\mathbf{x}^{k}}
\newcommand{\xkk}{\mathbf{x}^{k+1}}
\newcommand{\z}{\mathbf{z}}
\newcommand{\zk}{\mathbf{z}^{k}}
\newcommand{\zkk}{\mathbf{z}^{k+1}}

\renewcommand{\b}{\mathbf{b}}
\newcommand{\bflambda}{\mathbf{\lambda}}
\newcommand{\blambda}{\bm{\lambda}}
\newcommand{\blambdak}{\bm{\lambda}^{k}}
\newcommand{\blambdakk}{\bm{\lambda}^{k+1}}
\newcommand{\hatlambdakk}{\hat{\bm{\lambda}}^{k+1}}
\newcommand{\A}{\mathcal{A}}

\newcommand{\X}{\mathbf{X}}

\newcommand{\W}{\mathbf{W}}
\newcommand{\Z}{\mathbf{Z}}
\newcommand{\D}{\mathbf{D}}

\newcommand{\lbar}{\left\|}
\newcommand{\rbar}{\right\|}
\newcommand{\sumi}{\sum_{i=1}^n}

\DeclareMathOperator*{\argmin}{argmin}

\newcommand{\thetak}{{\theta}^{(k)}}
\newcommand{\thetakk}{{\theta}^{(k+1)}}
\newcommand{\betak}{{\beta}^{(k)}}
\newcommand{\betakk}{{\beta}^{(k+1)}}

\newtheorem{theorem}{Theorem}
\newtheorem{proposition}{Proposition}

\newtheorem{lemma}{Lemma}

\renewcommand{\a}{\mathbf{a}}
\newcommand{\ak}{a^{(k)}}
\newcommand{\akk}{a^{(k+1)}}
\newcommand{\bk}{b^{(k)}}
\newcommand{\bkk}{b^{(k+1)}}
\renewcommand{\b}{\mathbf{b}}
\newcommand{\sumk}{\sum_{k=0}^K}

\begin{document}
%
\title{Fast Proximal Linearized Alternating Direction Method of Multiplier with Parallel Splitting}
	
\author{Canyi Lu$^1$, Huan Li$^2$, Zhouchen Lin$^{2,3,}$\thanks{Corresponding author.}, Shuicheng Yan$^1$\\
	$^1$ Department of Electrical and Computer Engineering, National University of Singapore\\
	$^2$ Key Laboratory of Machine Perception (MOE), School of EECS, Peking University\\
	$^3$ Cooperative Medianet Innovation Center, Shanghai Jiaotong University \\
	{\tt\small canyilu@gmail.com, lihuanss@pku.edu.cn, zlin@pku.edu.cn, eleyans@nus.edu.sg}
}

\maketitle
\vspace{-100em}
\begin{abstract}
\begin{quote}
 The Augmented Lagragian Method (ALM) and Alternating Direction Method of Multiplier (ADMM) have been powerful optimization methods for general convex programming subject to linear constraint. We consider the convex problem whose objective consists of a smooth part and a nonsmooth but simple part. We propose the Fast Proximal Augmented Lagragian Method (Fast PALM) which achieves the convergence rate $O(1/K^2)$, compared with  $O(1/K)$ by the traditional PALM. In order to further reduce the per-iteration complexity and handle the multi-blocks problem, we propose the Fast Proximal ADMM with Parallel Splitting (Fast PL-ADMM-PS) method. It also partially improves the rate related to the smooth part of the objective function. Experimental results on both synthesized and real world data demonstrate that our fast methods significantly improve the previous PALM and ADMM. 
\end{quote}
\end{abstract}
\section{Introduction}

This work aims to solve the following linearly constrained separable convex problem with $n$ blocks of variables
\begin{equation}
\begin{split}
\min\limits_{\x_1,\cdots,\x_n} & \  f(\x)=\sum\limits_{i=1}^nf_i(\x_i)=\sumi\left( g_i(\x_i)+h_i(\x_i)\right), \\
\text{s.t.} & \quad
\A(\x)=\sum\limits_{i=1}^n \A_i(\x_i)=\b,\label{eq:model_problem_multivar}
\end{split}
\end{equation}
where $\x_i$'s and $\b$ can be vectors or matrices and both $g_i$ and $h_i$ are convex and lower semi-continuous. For $g_i$, we assume that $\nabla g_i$ is Lipschitz continuous with the Lipschitz constant $L_i>0$, i.e, $\lbar \nabla g_i(\x_i)- \nabla g_i(\y_i)\rbar \leq L_i \lbar \x_i-\y_i\rbar, \forall \x_i,\y_i$. For $h_i$, we assume that it may be nonsmooth and it is simple, in the sense that the proximal operator problem $\min_{\x} h_i(\x) + \frac{\alpha}{2}||\x - \a||^2$ ($\alpha>0$) can be cheaply solved. The bounded mappings $\A_i$'s are linear (e.g., linear transformation or the sub-sampling operator in matrix completion \cite{candes2009exact}). For the simplicity of discussion, we denote $\x = [\x_1 ; \x_2; \cdots ; \x_n ]$, $\A = [\A_1, \A_2, \cdots, \A_n]$ and
$\sum_{i=1}^{n} \A_i(\x_i) = \A(\x)$, $f_i=g_i+h_i$. For any compact set $X$, let $D_{{X}} =
\sup_{\x_1,\x_2 \in {X}} ||\x_1 - \x_2||$ be
the diameter of ${X}$. We also denote $\D_{\x^{*}} =
||\x^{0} - \x^{*}||$. We assume there exists a saddle point $(\x^{*}, \bm{\lambda}^{*}) \in X \times {\Lambda}$ to (\ref{eq:model_problem_multivar}), i.e., $\A(\x^{*}) = \b$ and $- \A_i^T(\bm{\lambda}^{*}) \in \partial f_i(\x_i^{*}), \  i = 1,\cdots,n$, where $\A^T$ is the adjoint
operator of $\A$, $X$ and $\Lambda$ are the feasible sets of the primal variables and dual variables, respectively.

By using different $g_i$'s and $h_i$'s, a variety of machine learning problems can be cast into (\ref{eq:model_problem_multivar}), including Lasso \cite{tibshirani1996regression} and its variants \cite{lu2013correlation,jacob2009group}, low rank matrix decomposition \cite{RPCA}, completion \cite{candes2009exact} and representation model \cite{lu2012robust,liu2011latent} and latent variable graphical model selection~\cite{chandrasekaran2010latent}. Specifically, examples of $g_i$ are: (i) the square loss $\frac{1}{2}||\D \x - \y||^2$ , where $\D$ and $\y$ are of compatible dimensions. A more special case is the known Laplacian regularizer $\text{Tr}(\X \mathbf{L} \X^T)$, where $\mathbf{L}$ is the Laplacian matrix which is positive semi-definite; (ii) Logistic loss $\sum^{m}_{i=1} \log (1 + \exp(-y_i \mathbf{d}_i^T \x))$, where $\mathbf{d}_i$'s and $y_i$'s are the data points and the corresponding labels, respectively; (iii) smooth-zero-one loss $\sum^{m}_{i=1} \frac{1}{1+ \exp(c y_i \mathbf{d}_i^T \x)}$, $c > 0$.
The possibly nonsmooth $h_i$ can be many norms, e.g., $\ell_1$-norm $||\cdot||_{1}$ (the sum of absolute values of all entries), $\ell_2$-norm $|| \cdot ||$ or Frobenius norm $||\cdot||_F$ and nuclear norm $|| \cdot ||_{*}$ (the sum of the singular values of a matrix).

This paper focuses on the popular approaches which study problem (\ref{eq:model_problem_multivar}) from the aspect of the  augmented Lagrangian function $L(\x,\bm{\lambda})  = f(\x) + \langle \bm{\lambda}, \A(\x) - \b \rangle + \frac{\beta}{2} ||\A(\x) - \b ||^2$, where $\bm{\lambda}$ is the Lagrangian multiplier or dual variable and $\beta > 0$. A basic idea to solve problem (\ref{eq:model_problem_multivar}) based on $L(\x,\bm{\lambda})$ is the Augmented Lagrangian Method (ALM) \cite{hestenes1969multiplier}, which is a special case of the Douglas-Rachford splitting \cite{douglas1956numerical}. 

\begin{table*}[t]
	\begin{tabular}{c|c|c|c}
		\hline
		{ PALM } & { Fast PALM} & { PL-ADMM-PS} & { Fast PL-ADMM-PS}  \\
		\hline
		\begin{minipage}{2.5cm} \centering \vspace{0.2cm}
			$O\left(\frac{ D^2_{\x^{*}} + D^2_{\bm{\lambda}^{*}} }{K}  \right)$ 
			\vspace{0.1cm}\end{minipage} &
		\begin{minipage}{2.5cm} \centering \vspace{0.2cm} $O\left(\frac{ D^2_{\x^{*}} + D^2_{\bm{\lambda}^{*}} }{K^2}  \right)$
			\vspace{0.1cm}\end{minipage} & \begin{minipage}{3.5cm} \centering \vspace{0.2cm} $O\left( \frac{D^2_{{\x}^*}}{K} + \frac{D^2_{\x^*}}{K} + \frac{D^2_{\bm{\lambda}^*}}{ K} \right)$ 
			\vspace{0.1cm}\end{minipage}  & 
		\begin{minipage}{3.5cm}  \centering \vspace{0.2cm} $O \left( \frac{ D^2_{{\x}^*}}{K^2} + \frac{D^2_{{X}}}{K} + \frac{D^2_{\Lambda}}{ K} \right)$  \vspace{0.1cm}\end{minipage}  \\
		\hline
	\end{tabular}	
	\centering \caption{Comparison of the convergence rates of previous methods and our fast versions}	\label{all_convergence_rates}
	\vspace{-0.6em}
\end{table*}

An influential variant of ALM is the Alternating Direction Mehtod of Multiplier (ADMM) \cite{boyd2011distributed}, which solves
problem (\ref{eq:model_problem_multivar})
with $n = 2$ blocks of variables.
However, the cost for solving the subproblems in ALM and ADMM in each iteration is usually high when $f_i$ is not simple and  $\A_i$ is non-unitary ($\A_i^T\A_i$ is not the identity mapping). To alleviate this issue, the Linearized ALM (LALM) \cite{yang2013linearized} and Linearized ADMM (LADMM) \cite{LADMAP} were proposed by linearizing the augmented term $\frac{\beta}{2}||\A (\x) - \b||^2 $ and thus the subproblems are easier to solve. For (\ref{eq:model_problem_multivar}) with $n>2$ blocks of variables, the Proximal Jacobian ADMM \cite{Mintao} and Linearized ADMM with Parallel Splitting (L-ADMM-PS) \cite{LADMPS} 
guaranteed to solve (\ref{eq:model_problem_multivar}) when $g_i = 0$ with convergence guarantee. To further exploit the Lipschitz continuous gradient property of $g_i$'s in (\ref{eq:model_problem_multivar}), the 
work \cite{LADMPS} proposed a Proximal Linearized ADMM with Parallel Splitting (PL-ADMM-PS) by further linearizing the smooth part $g_i$. PL-ADMM-PS requires lower per-iteration cost than L-ADMM-PS for solving the general problem (\ref{eq:model_problem_multivar}).

%


Beyond the per-iteration cost, another important way to measure the speed of the algorithms is the convergence rate. Several previous work proved the convergence rates of the augmented Lagrangian function based methods \cite{he20121,Mintao,LADMPS}. Though the convergence functions used to measure the convergence rate are different, the 
convergence rates of all the above discussed methods for (\ref{eq:model_problem_multivar}) are all $O({1}/{K})$, where $K\textsc{}$ is the number of iterations.
However, the rate $O(1/K)$ may be suboptimal in some cases. Motivated by the seminal work \cite{nesterov1983method}, several fast first-order
methods with the optimal rate $O(1/K^2)$ have been developed for unconstrained 
problems \cite{beck2009fast,tseng}. More recently, by applying a similar accelerating technique, several fast ADMMs have been proposed to solve a special case of problem (\ref{eq:model_problem_multivar}) with $n=2$ blocks of variables
\begin{equation}\label{pro2222}
\min_{\x_1,\x_2} g_1(\x_1) + h_2(\x_2),\ \
\text{s.t.}\ \ \A_1 (\x_1) + \A_2 (\x_2) = \b.
\end{equation}
A fast ADMM proposed in \cite{azadi2014towards}\footnote{The method in \cite{azadi2014towards} is a fast stochastic ADMM. It is easy to give the corresponding deterministic version by computing the gradient in each iteration exactly to solve (\ref{pro2222}).} is able to solve (\ref{pro2222}) with the convergence rate $O\left(\frac{D_{X}^2}{K^2} + \frac{D_\Lambda^2}{K}\right)$. But their
result is a bit weak since their used function to characterize the
convergence can be negative. The work \cite{ouyang2014accelerated} proposed another
fast ADMM with the rate $O\left(\frac{ D_{\x^{*}}^2}{K^2} + \frac{ D_{\x^{*}}^2}{K}\right)$
for primal residual and $O\left(\frac{ D_{\x^{*}}^2}{K^{3/2}} + \frac{D_{\x^{*}} +  D_{\bm{\lambda}^{*}}}{K}\right)$ for feasibility
residual. However, their result requires that the number of iterations $K$ should be predefined, which is not reasonable in practice. It is usually difficult in practice to determine the optimal $K$ since we usually stop the algorithms when both the primal and feasibility residuals are sufficiently small \cite{LADMPS}. The fast ALM proposed in \cite{he2010acceleration} owns the convergence rate $O(1/K^2)$, but it requires the objective $f$ to be differentiable. This limits its applications for nonsmooth optimization in most compressed sensing problems. Another work \cite{goldstein2012fast} proved a better convergence rate than $O(1/K)$ for ADMM. But their method requires much stronger assumptions, e.g., strongly convexity of $f_i$'s, which are usually violated in practice. In this work, we only consider  (\ref{eq:model_problem_multivar}) whose objective is not necessarily strongly convex.

In this work, we aim to propose fast ALM type methods to solve  the general problem (\ref{eq:model_problem_multivar}) with optimal convergence rates. The contributions are summarized as follows:
\begin{itemize}
	\item First, we consider (\ref{eq:model_problem_multivar}) with $n=1$ (or one may regard all $n$ blocks as a superblock) and propose the Fast Proximal Augmented Lagrangian Method (Fast PALM). We prove that Fast PALM converges with the rate
	$O\left(\frac{ D^2_{\x^{*}} +  D^2_{\bm{\lambda}^{*}}}{K^2}\right)$, which is a significant improvement of ALM/PALM\footnote{PALM is a variant of ALM proposed in this work.} with rate $O\left(\frac{  D^2_{\x^{*}} +  D^2_{\bm{\lambda}^{*}}}{K}\right)$. To the best of our knowledge, Fast PALM is the first improved ALM/PALM which achieves the  rate $O(1/K^2)$ for the nonsmooth problem (\ref{eq:model_problem_multivar}).
	\item Second, we consider (\ref{eq:model_problem_multivar}) with $n>2$ and propose the Fast Proximal Linearized
	ADMM with Parallel Splitting (Fast PL-ADMM-PS), which converges with rate $O \left( \frac{ D^2_{{\x}^*}}{K^2} + \frac{D^2_{{X}}}{K} + \frac{D^2_{\Lambda}}{ K} \right)$. As  discussed in Section 1.3 of \cite{ouyang2014accelerated}, such a rate is optimal and thus  is better than PL-ADMM-PS with rate $O\left( \frac{D^2_{{\x}^*}}{K} + \frac{D^2_{\x^*}}{K} + \frac{D^2_{\bm{\lambda}^*}}{ K} \right)$   \cite{LADMPS}. To the best of our knowledge, Fast PL-ADMM-PS is the first fast Jacobian type (update the variables in parallel) method to solve (\ref{eq:model_problem_multivar}) when $n>2$ with convergence guarantee. 
\end{itemize}
Table \ref{all_convergence_rates} shows the comparison of the convergence rates of previous methods and our fast versions. Note that Fast PALM and Fast PL-ADMM-PS have the same pter-iteration cost as PALM and PL-ADMM-PS, respectively. But the per-iteration cost of PL-ADMM-PS and Fast PL-ADMM-PS may be much cheaper than PALM and Fast PALM.


\section{Fast Proximal Augmented Lagrangian Method}

In this section, we consider (\ref{eq:model_problem_multivar}) with $n=1$ block of variable,
\begin{equation}\label{pro_one}
\min_{\x} f(\x) = g(\x) + h(\x),\quad
\text{s.t.}\quad \A (\x) = \b,
\end{equation}
where $g$ and $h$ are convex and $\nabla g$ is Lipschitz continuous with the Lipschitz constant $L$. The above problem can be solved by the traditional ALM which updates $\x$ and $\bm{\lambda}$ by
\begin{equation}
\left\{
\begin{aligned}
\x^{k+1}=&\arg\min_{\x} g(\x)+h(\x)+\langle \bm{\lambda}^k, \A(\x) - \b \rangle \\
& +\frac{\betak}{2} ||\A(\x) - \b ||^2 ,\label{updatexalm}\\
\bm{\lambda}^{k+1} =& \bm{\lambda}^k + \betak (\A(\x^{k+1}) - \b),
\end{aligned}
\right.
\end{equation}
where $\betak>0$. Note that $\nabla g$ is Lipschitz continuous. We have \cite{nesterov2004introductory}
\begin{equation}\label{propg}
g(\x)\leq g(\x^k)+\langle \nabla g(\x^k), \x - \x^k \rangle+\frac{L}{2} ||\x - \x^k||^2.
\end{equation}
This motivates us to use the right hand side of (\ref{propg}) as a surrogate of $g$ in (\ref{updatexalm}). Thus we can update $\x$ by solving the following problem which is simpler than (\ref{propg}),
\begin{equation*}
\begin{split}
&\x^{k+1}= \arg\min_{\x} g(\x^k)+\langle \nabla g(\x^k), \x - \x^k \rangle + h(\x) \\
&+ \langle \bm{\lambda}^k, \A(\x) - \b \rangle +
\frac{\betak}{2} ||\A(\x) - \b ||^2+ \frac{L}{2} ||\x - \x^k||^2. 
\end{split}
\end{equation*}
We call the method by using the above updating rule as Proximal Augmented Lagrangian Method (PALM). PALM can be regarded as a special case of Proximal Linearized Alternating Direction Method of Multiplier with Parallel Splitting in \cite{LADMPS} and it owns the convergence rate $O\left({1}/{K}\right)$, which is the same as the traditional ALM and ADMM. However, such a rate is suboptimal. Motivated by 
the technique from the accelerated proximal gradient method \cite{tseng}, we propose the Fast PALM as shown in Algorithm~\ref{alg1}. It uses the interpolatory sequences $\y^k$ and $\z^k$ as well as the stepsize
$\thetak$. Note that if we set $\thetak = 1$ in each iteration, Algorithm \ref{alg1} reduces to PALM. With careful choices of $\thetak$ and 
$\betak$ in Algorithm \ref{alg1}, we can accelerate the convergence rate of PALM from $O\left({1}/{K}\right)$ to $O({1}/{K^2})$.

\begin{algorithm}[t]\label{alg1}
	\caption{Fast PALM Algorithm}
	\hrule
	\hrule
	\vspace{0.1cm}
	\textbf{Initialize}: $\x^0$, $\z^0$, $\bm{\lambda}^0$, $\beta^{(0)}=\theta^{(0)}=1$.\\
	\For{$k = 0, 1, 2,\cdots$}{
		\begin{align}
		\ykk&=(1-\thetak)\x^k+\thetak\zk;\\
		\z^{k+1}&=\argmin\limits_{\x} \left \langle\nabla g(\y^{k+1}),\x\right \rangle+h(\x)\notag\\
		&+ \left \langle\bm\lambda^k,\A(\x) \right \rangle+\frac{\betak}{2}\|\A(\x)-\b\|^2  \notag\\
		&+\frac{L\thetak}{2}\|\x-\z^k\|^2;\label{updatexkfpalm} \\
		\x^{k+1}&=(1-\thetak)\x^k+\thetak\z^{k+1}; \label{ddefinex}\\
		\bm{\lambda}^{k+1}&=\bm{\lambda}^k+\betak(\A(\z^{k+1})-\b);\label{update_lambda_one} \\
		\thetakk&=\frac{-(\thetak)^2+\sqrt{(\thetak)^4+4(\thetak)^2}}{2}; \\
		\betakk&= \frac{1}{\thetakk}.
		\end{align}
	}
	\hrule
	\hrule\hrule
	\vspace{0.1cm}
\end{algorithm}

\begin{proposition}\label{prop_one}
In Algorithm \ref{alg1},  for any $\x$, we have
	\begin{align}
	&\frac{1-\thetakk}{(\thetakk)^2} \left(f(\x^{k+1})-f(\x)\right)-\frac{1}{\thetak}\left \langle\A^T(\blambdakk),\x-\z^{k+1} \right \rangle\qquad \notag\\
	&\leq\frac{1-\thetak}{(\thetak)^2}\left(f(\x^k)-f(\x)\right) \label{eqpro1}\\
	& \ \ \  \ + \frac{L}{2}\left(\|\z^{k}-\x\|^2-\|\z^{k+1}-\x\|^2\right).\notag
	\end{align}
\end{proposition}

\begin{theorem}\label{conv_rate_one}
	In Algorithm \ref{alg1}, for any $K>0$, we have
	\begin{align}\label{eqconrateone}
	&	f(\x^{K+1})-f(\x^*)+\left \langle\bm\lambda^*,\A(\x^{K+1})-\b \right \rangle+\frac{1}{2}\|\A (\x^{K+1})-\b\|^2 \notag\\
	&\leq\frac{2}{(K+2)^2}\left(LD_{\x^*}^2+D_{\blambda^*}^2\right).
	\end{align}
\end{theorem}
We use the convergence function, i.e., the left hand side of (\ref{eqconrateone}), in \cite{LADMPS} to measure the convergence rate of the algorithms in this work. Theorem \ref{conv_rate_one} shows that our Fast PALM achieves the rate $O\left(\frac{LD_{\x^*}^2+D^2_{\blambda^*}}{K^2}\right)$, which is much better than $O\left(\frac{LD^2_{\x^*}+\frac{1}{\beta}D^2_{\blambda^*}}{K}\right)$ by PALM\footnote{It is easy to achieve this since PALM is a special case of Fast PALM by taking $\thetak=1$. }.
The improvement of Fast PALM over PALM is similar to the one of Fast ISTA over ISTA \cite{beck2009fast,tseng}. The difference is that Fast ISTA targets for unconstrained problem which is easier than our problem (\ref{eq:model_problem_multivar}). Actually, if the constraint in (\ref{eq:model_problem_multivar}) is dropped (i.e., $\A=\mathbf{0}$, $\b=\mathbf{0}$), our Fast PALM is similar as the Fast ISTA. 

We would like to emphasize some key differences between our Fast PALM and previous fast ALM type methods \cite{azadi2014towards,ouyang2014accelerated,he2010acceleration}. First, it is easy to apply the two blocks  fast ADMM methods in \cite{azadi2014towards,ouyang2014accelerated} to solve problem (\ref{pro_one}). Following their choices of parameters and proofs, the convergence rates are still $O({1}/{K})$. The key improvement of our method comes from the different choices of 
$\thetak$ and $\betak$ as shown in Theorem \ref{conv_rate_one}. The readers can refer to the detailed proofs in the supplementary material. Second, the fast ADMM in \cite{ouyang2014accelerated} requires  predefining 
the total number of iterations, which is usually difficult
in practice. However, our Fast PALM has no such a limitation. Third, the fast ALM in \cite{he2010acceleration} also owns the rate $O(1/K^2)$. But it is restricted to differentiable objective minimization and thus is not applicable to our problem (\ref{eq:model_problem_multivar}). Our method has no such a limitation.  

A main limitation of PALM and Fast PALM is that their per-iteration cost may be high when $h_i$ is nonsmooth and $\A_i$ is non-unitary. In this case, solving the subproblem (\ref{updatexkfpalm}) requires calling other iterative solver, e.g., Fast ISTA \cite{beck2009fast}, and thus the high per-iteration cost may limit the application of Fast PALM. In next section, we present a fast ADMM which has lower per-iteration cost.

\section{Fast Proximal Linearized  ADMM with Parallel Splitting}
In this section, we consider problem (\ref{eq:model_problem_multivar}) with $n>2$ blocks of variables. The state-of-the-art solver for (\ref{eq:model_problem_multivar}) is the Proximal Linearized ADMM with Parallel Splitting (PL-ADMM-PS) \cite{LADMPS} which updates each $\x_i$ in parallel by
\begin{align}
\x_i^{k+1}&=\argmin\limits_{\x_i} g_i(\x_i^{k})+ \left \langle\nabla g_i(\x_i^{k}),\x_i-\x_i^{k} \right \rangle+h_i(\x_i) \notag \\
&+ \left \langle\bm{\lambda}^k,\A_i(\x_i) \right \rangle+\left\langle \betak\A_i^T\left( \A(\x^k)-\b\right),\x_i-\x_i^k\right\rangle \notag\\
&+\frac{L_i+\betak\eta_i}{2}\|\x_i-\x_i^k\|^2,\label{updatexladmmpsss2}
\end{align}
where $\eta_i>n||\A_i||^2$ and $\betak>0$. Note that the subproblem (\ref{updatexladmmpsss2}) is easy to solve when $h_i$ is nonsmooth but simple. Thus PL-ADMM-PS has much lower per-iteration cost than PALM and Fast PALM. On the other hand, PL-ADMM-PS converges with the rate $O(1/K)$ \cite{LADMPS}. However, such a rate is also suboptimal. Now we show that it can be further accelerated by a similar technique as that in Fast PALM. 
See Algorithm \ref{alg2} for our Fast PL-ADMM-PS.

\begin{proposition}\label{prop3}
	In Algorithm \ref{alg2}, for any $\x_i$, we have
	\begin{align}
	&\frac{1-\thetakk}{(\thetakk)^2}\left(f_i(\x_i^{k+1})-f_i(\x_i)\right) \notag\\
	& -\frac{1}{\thetak} \left \langle\A_i^T({\hatlambdakk}),\x_i-\z_i^{k+1} \right \rangle\notag\\
	\leq&\frac{1-\thetak}{(\thetak)^2}\left(f_i(\x_i^k)-f_i(\x_i)\right) \label{pro331} \\
	& +\frac{L_i}{2}\left(\|\z_i^{k}-\x_i\|^2-\|\z_i^{k+1}-\x_i\|^2\right)\notag\\
	&+\frac{\betak\eta_i}{2\thetak}\left(\|\z_i^{k}-\x_i\|^2-\|\z_i^{k+1}-\x_i\|^2-\|\z_i^{k+1}-\z_i^k\|^2\right),\notag
	\end{align}
	where
	$\hat{\bm\bflambda}^{k+1}=\bm\lambda^k+\betak\left(\A(\z^{k})-\b\right)$.
\end{proposition}
\begin{theorem}\label{con_fastPS}
	In Algorithm \ref{alg2}, for any $K>0$, we have
	\begin{align}
	&  f(\x^{K+1})-f(\x^*)+\left \langle\bm{\lambda}^*,\A(\x^{K+1})-\b \right \rangle\notag\\
	& +\frac{\beta\alpha}{2}\left\|\A(\x^{K+1})-\b\right\|^2 \label{con_FastPS}\\
	\leq  & \frac{2L_{\max}D^2_{\x^*}}{(K+2)^2}+\frac{2\beta\eta_{\max} D^2_{X}}{K+2}+\frac{2D^2_{\Lambda}}{\beta(K+2)},\notag
	\end{align}
	where $\alpha=\min\left\{\frac{1}{n+1},\left\{\frac{\eta_i-n\|\A_i\|^2}{2(n+1)\|\A_i\|^2},i=1,\cdots,n\right\}\right\}$,  $L_{\max}=\max\{L_i,i=1,\cdots,n\}$ and $\eta_{\max}=\max\{\eta_i,i=1,\cdots,n\}$.
\end{theorem}

\newcommand{\parallelsum}{\mathbin{\!/\mkern-5mu/\!}}
\begin{algorithm}[!t]\label{alg2}
	\DontPrintSemicolon
	\caption{Fast PL-ADMM-PS Algorithm}
	\hrule
	\hrule
	\vspace{0.1cm}
	\textbf{Initialize}: $\x^0$, $\z^0$, $\bm{\lambda}^0$, $\theta^{(0)}=1$, fix $\betak=\beta$ for  $k\geq0$, $\eta_i>n\|\A_i\|^2$, $i=1,\cdots,n$,\\
	\For{$k = 0, 1, 2, \cdots$}{
		\hspace*{0.1cm}$\parallelsum$ Update $\y_i$, $\z_i$, $\x_i$, $i=1,\cdots,n$, in parallel by
\begin{align}
\y_i^{k+1}&=(1-\thetak)\x_i^k+\thetak\z_i^k;\\
\z_i^{k+1}&=\argmin\limits_{\x_i}  \left \langle\nabla g_i(\y_i^{k+1}),\x_i \right \rangle +h_i(\x_i)\notag \\
&+ \left \langle\bm{\lambda}^k,\A_i(\x_i) \right \rangle+\left\langle \betak\A_i^T\left( \A(\z^k)-\b\right),\x_i\right\rangle   \notag \\
&+\frac{L(g_i)\thetak+\betak\eta_i}{2}\|\x_i-\z_i^k\|^2; \label{update_z}\\
\x_i^{k+1}&=(1-\thetak)\x_i^k+\thetak\z_i^{k+1}; \\
\mathbf{\bm{\lambda}}^{k+1}&=\mathbf{\bm{\lambda}}^k+\beta^k\left(\A(\z^{k+1})-\b\right);\label{update_lambda}\\
\thetakk&=\frac{-(\thetak)^2+\sqrt{(\thetak)^4+4(\thetak)^2}}{2}.
\end{align}
	}	
	\hrule
	\hrule
	\hrule
	\vspace{0.1cm}
\end{algorithm}

%
From Theorem \ref{con_fastPS}, it can be seen that our Fast PL-ADMM-PS \emph{partially} accelerates the convergence rate of PL-ADMM-PS from  $O\left(\frac{L_{\max}D_{\x^{*}}^2}{K} + \frac{\beta \eta_{\max}D^2_{\x^{*}}}{K} + \frac{D^2_{\bm{\lambda^{*}}}}{\beta K} \right)$ to $O\left(\frac{L_{\max} D_{\x^{*}}^2}{K^2} + \frac{\beta\eta_{\max}  D^2_{{X}}}{K} + \frac{D^2_{\Lambda}}{\beta K} \right)$. 
Although the improved rate is also $O(1/K)$, what makes it more attractive is that it allows very large Lipschitz constants $L_i$'s. In particular, $L_i$ can be as large as $O(K)$, without affecting the rate of convergence (up to a constant factor). 
The above improvement is the same as fast ADMMs \cite{ouyang2014accelerated} for problem (\ref{pro2222}) with only $n=2$ blocks. But it is inferior to the Fast PALM over PALM. The key difference is that Fast PL-ADMM-PS further linearizes the augmented term $\frac{1}{2}||\A(\x)-\b||^2$. 
This improves the efficiency for solving the subproblem, but slows down the convergence. Actually, when linearizing the augmented term, we have a new term with the factor $\betak\eta_i/\thetak$ in (\ref{pro331}) (compared with (\ref{eqpro1}) in Fast PALM). Thus (\ref{con_FastPS}) has a new term by comparing with that in (\ref{eqconrateone}). This makes the choice of $\betak$ in Fast PL-ADMM-PS different from the one in Fast PALM. Intuitively, it can be seen that a larger value of $\betak$ will increase the second terms of (\ref{con_FastPS}) and decrease the third term of (\ref{con_FastPS}). Thus $\betak$ should be fixed in order to guarantee the convergence. This is different from the choice of $\betak$ in Fast PALM which is adaptive to the choice of $\thetak$. 

Compared with PL-ADMM-PS, our Fast PL-ADMM-PS achieves a better rate, but with the price on the boundedness of the feasible primal set 
${X}$ and the feasible dual set $\Lambda$. Note that many previous work, e.g.,  \cite{he20121,azadi2014towards}, also require such a boundedness assumption when proving the convergence of ADMMs. In the following, we give some conditions which guarantee such a boundedness assumption.

\begin{figure*}[t]
	\setcounter{subfigure}{0}
	\subfigure[$m=100$, $n=300$]{
		\includegraphics[width=0.24\textwidth]{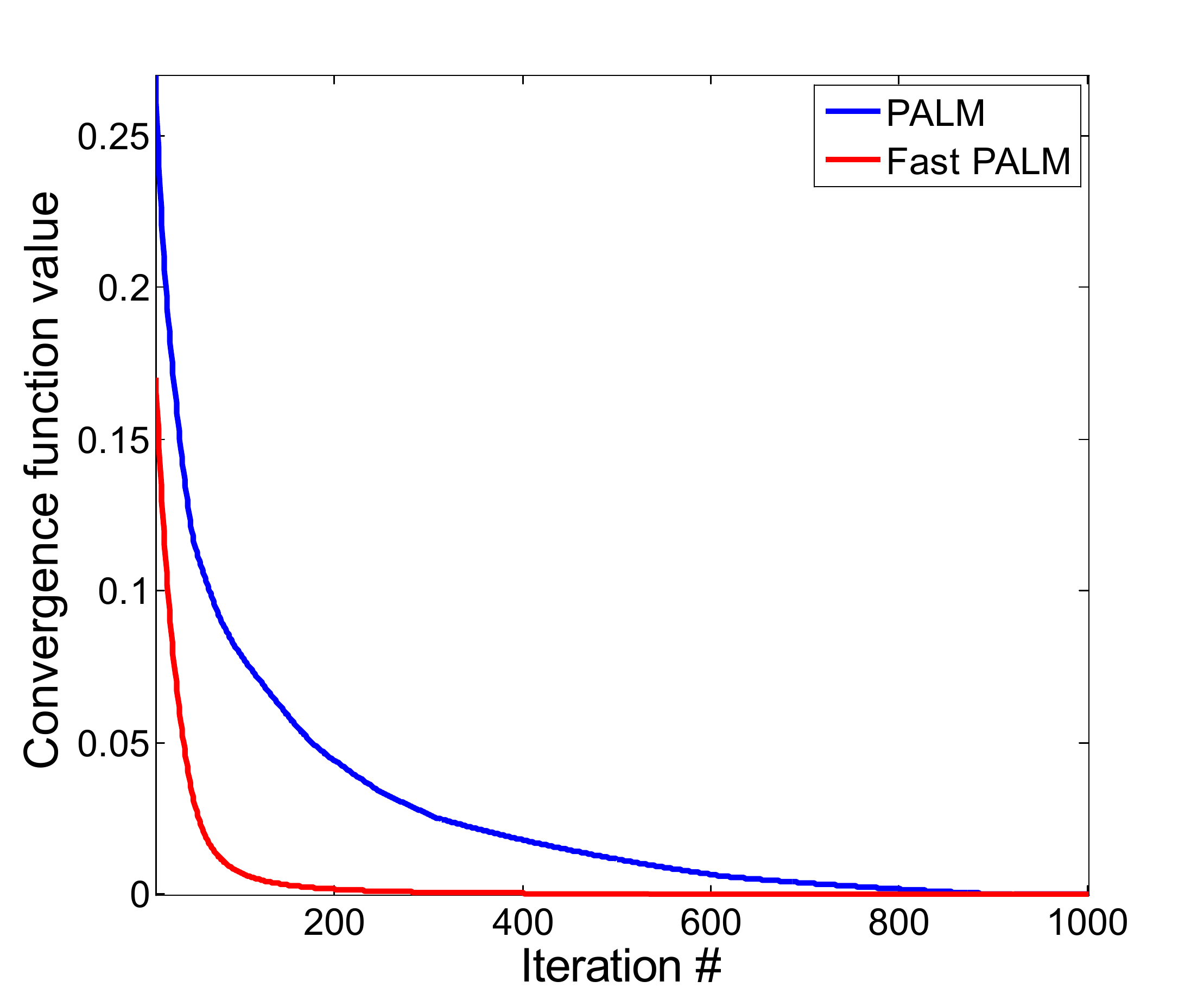}}
	\subfigure[$m=300$, $n=500$]{
		\includegraphics[width=0.24\textwidth]{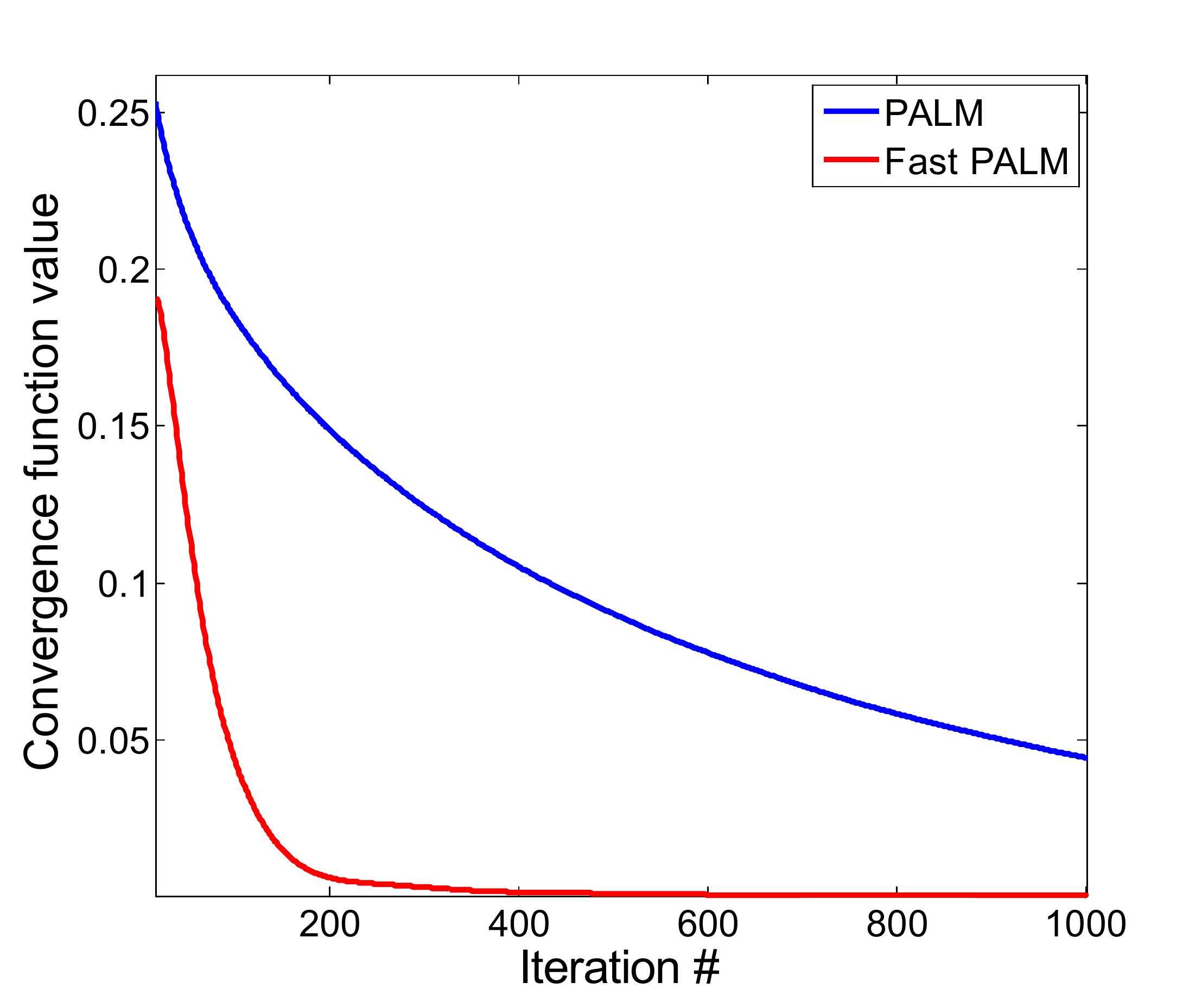}}
	\subfigure[$m=500$, $n=800$]{
		\includegraphics[width=0.24\textwidth]{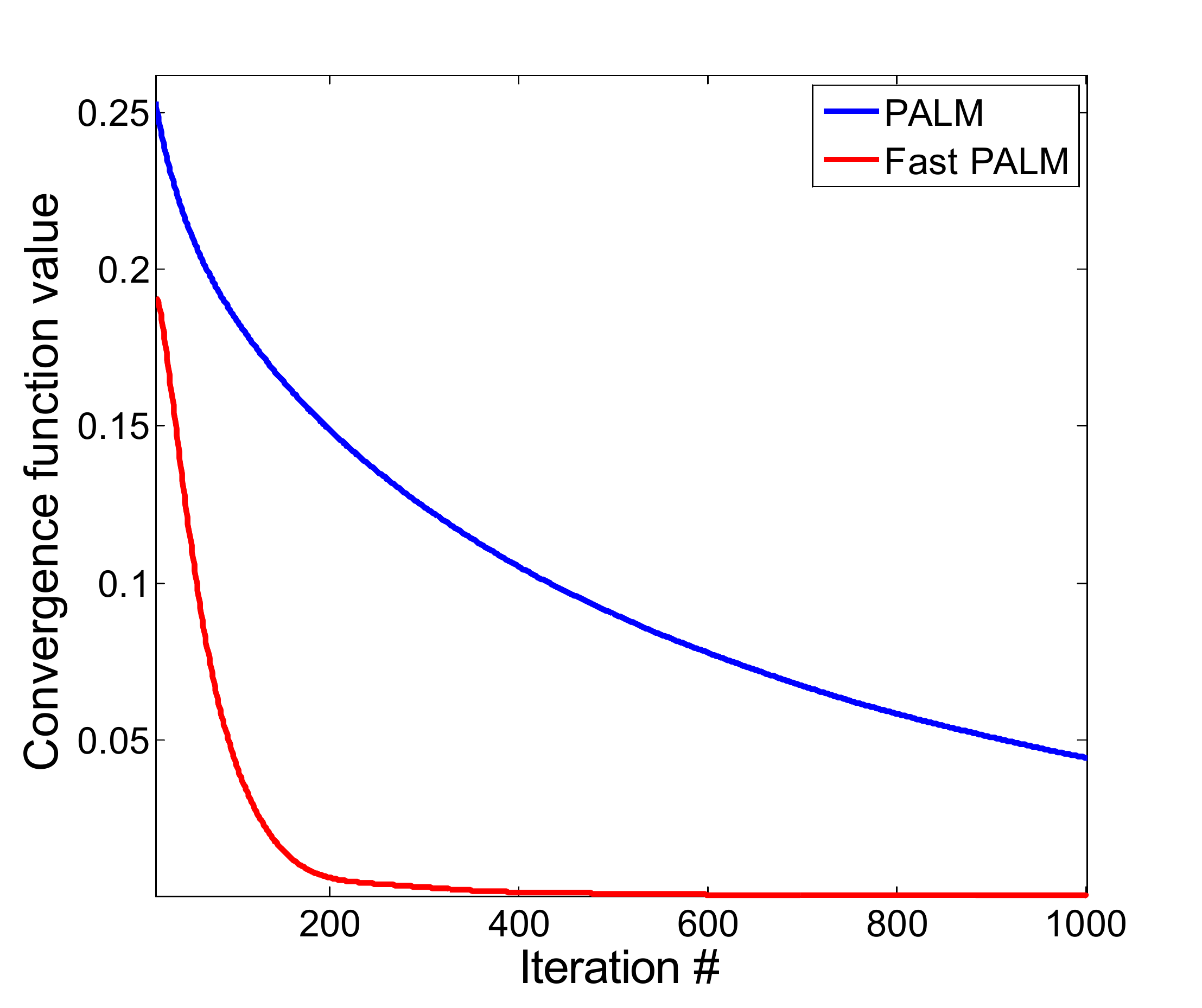}}
	\subfigure[$m=800$, $n=1000$]{
		\includegraphics[width=0.24\textwidth]{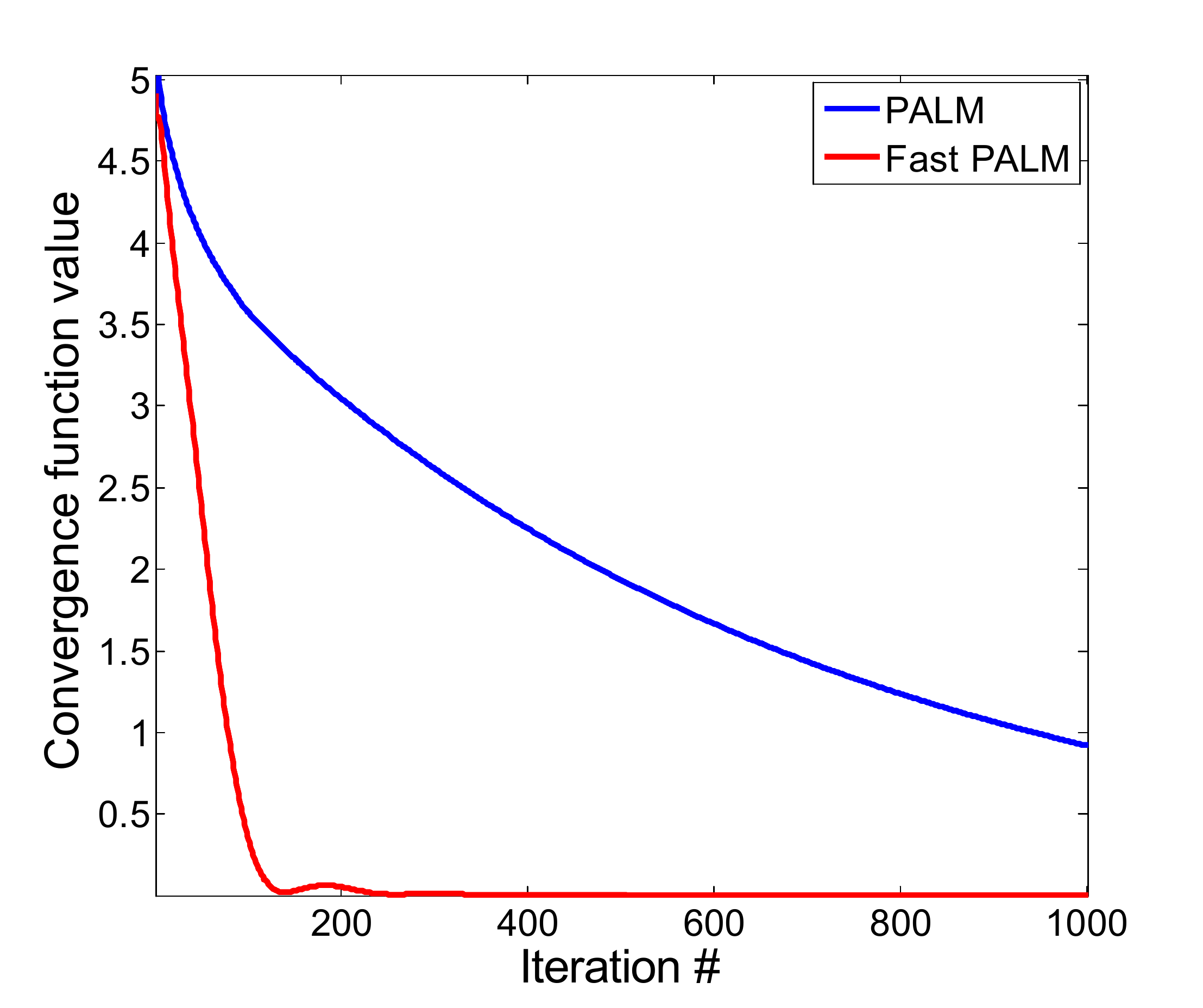}}\\
	\vspace{-0.3cm}
	\caption{Plots of the convergence function values of (\ref{eqconrateone}) in each iterations by using PALM and Fast PALM for (\ref{probl1}) with different sizes of $\mathbf{A}\in\mathbb{R}^{m\times n}$.} 
	\label{fig_res_fastALM}
	\vspace{-0.4cm}
\end{figure*}
\begin{theorem}\label{bound_lambda}
	Assume the mapping $\A(\x_1,\cdots,\x_n)=\sum_{i=1}^n\A_i(\x_i)$ is onto\footnote{This assumption is equivalent to that the matrix $A\equiv(A_1,\cdots,A_n)$ is of full row rank, where $A_i$ is the matrix representation of $\A_i$.}, the sequence $\{\z^k\}$ is bounded, $\partial h(\x)$ and $\nabla g(\x)$ are bounded if $\x$ is bounded, then $\{\x^k\}$, $\{\y^k\}$ and $\{\blambda^k\}$ are bounded.
\end{theorem}
Many convex functions, e.g., the $\ell_1$-norm, in compressed sensing own the bounded subgradient.

\section{Experiments}
In this section, we report some numerical results to demonstrate the effectiveness of our fast PALM and PL-ADMM-PS. We first compare our Fast PALM which owns the optimal convergence rate $O(1/K^2)$ with the basic PALM on a problem with only one block of variable. Then we conduct two experiments to compare our Fast PL-ADMM-PS with PL-ADMM-PS on two multi-blocks problems. The first one is tested on the synthesized data, while the second one is for subspace clustering tested on  the real-world data. We examine the convergence behaviors of the compared methods based on the convergence functions shown in (\ref{eqconrateone}) and (\ref{con_FastPS}).  All the numerical experiments are run on a PC with 8 GB of RAM and Intel Core 2 Quad CPU Q9550.
\subsection{Comparison of PALM and Fast PALM}
We consider the following problem
\begin{equation}\label{probl1}
\min_{\x} ||\x||_1+\frac{\alpha}{2}||\mathbf{A}\x-\b||_2^2, \ \ \text{s.t.} \ \ \mathbf{1}^T\x=1,
\end{equation}
where $\alpha>0$, $\mathbf{A}\in\mathbb{R}^{m\times n}$, $\b\in\mathbb{R}^m$, and $\mathbf{1}\in\mathbb{R}^n$ is the all one vector. There may have many fast solvers for problem (\ref{probl1}). In this experiment, we focus on the performance comparison of PALM and Fast PALM for (\ref{probl1}). Note that the per-iteration cost of these two methods are the same. Both of them requires solving an $\ell_1$-minimization problem in each iteration. In this work, we use the SPAMS package \cite{mairal2010online} to solve it which is very fast.

The data matrix $\mathbf{A}\in\mathbb{R}^{m\times n}$, and $\mathbf{b}\in\mathbb{R}^m$ are generated by the Matlab command \mcode{randn}. We conduct four experiments on different sizes of $\mathbf{A}$ and $\mathbf{b}$. We use the left hand side of (\ref{eqconrateone}) as the convergence function to evaluate the convergence behaviors of PALM and Fast PALM. For the saddle point $(\x^*,\blambda^*)$ in (\ref{eqconrateone}), we run the Fast PALM with 10,000 iterations and use the obtained solution as the saddle point. Figure \ref{fig_res_fastALM} plots the convergence functions value within 1,000 iterations. It can be seen that our Fast PALM  converges much faster than PALM. Such a result verifies our theoretical improvement of Fast PALM with optimal rate $O(1/K^2)$ over PALM with the rate $O(1/K)$.

\begin{figure*}[t]
	\setcounter{subfigure}{0}
	\subfigure[$m=100$]{
		\includegraphics[width=0.32\textwidth]{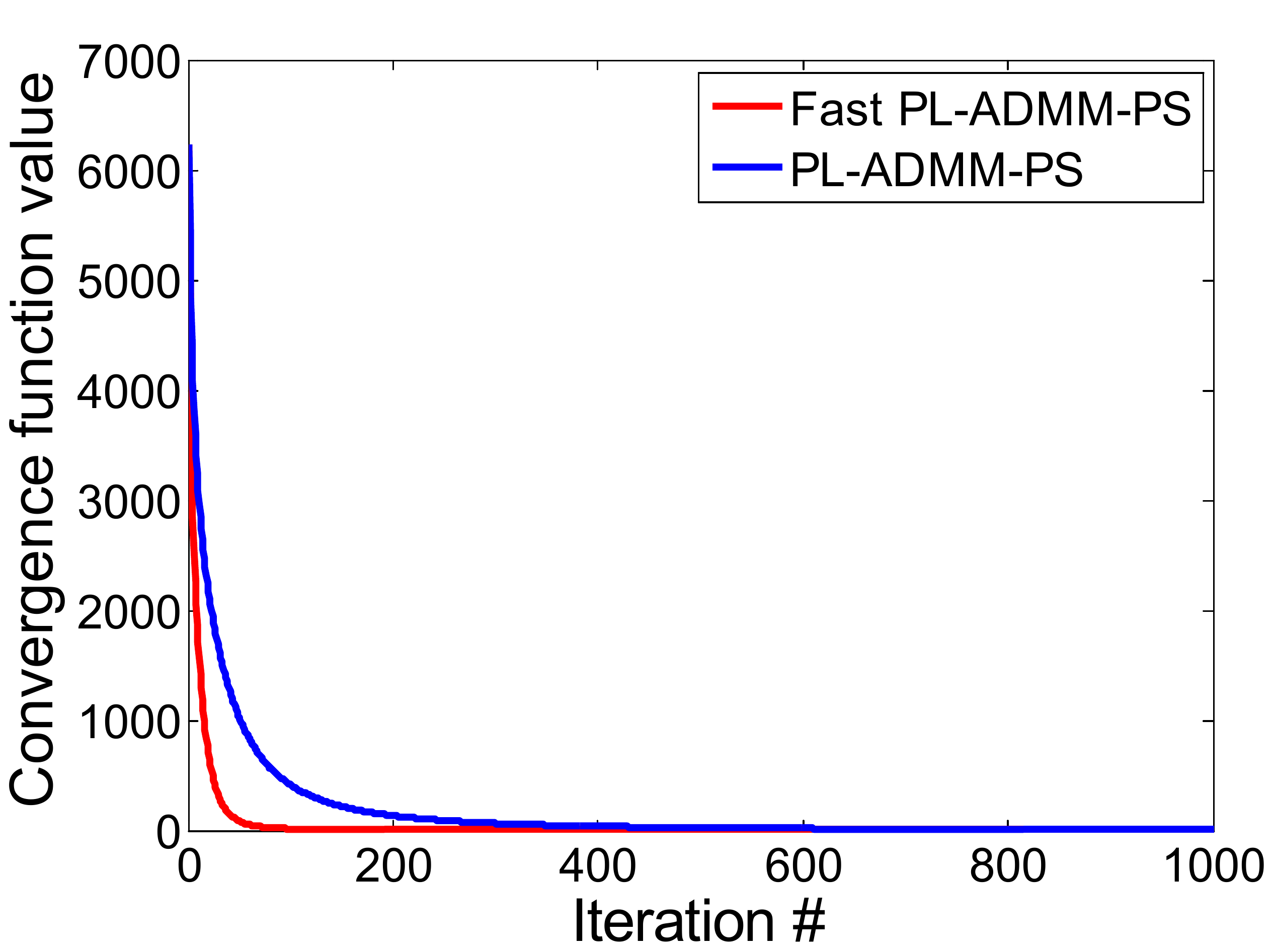}}
	\subfigure[$m=300$]{
		\includegraphics[width=0.31\textwidth]{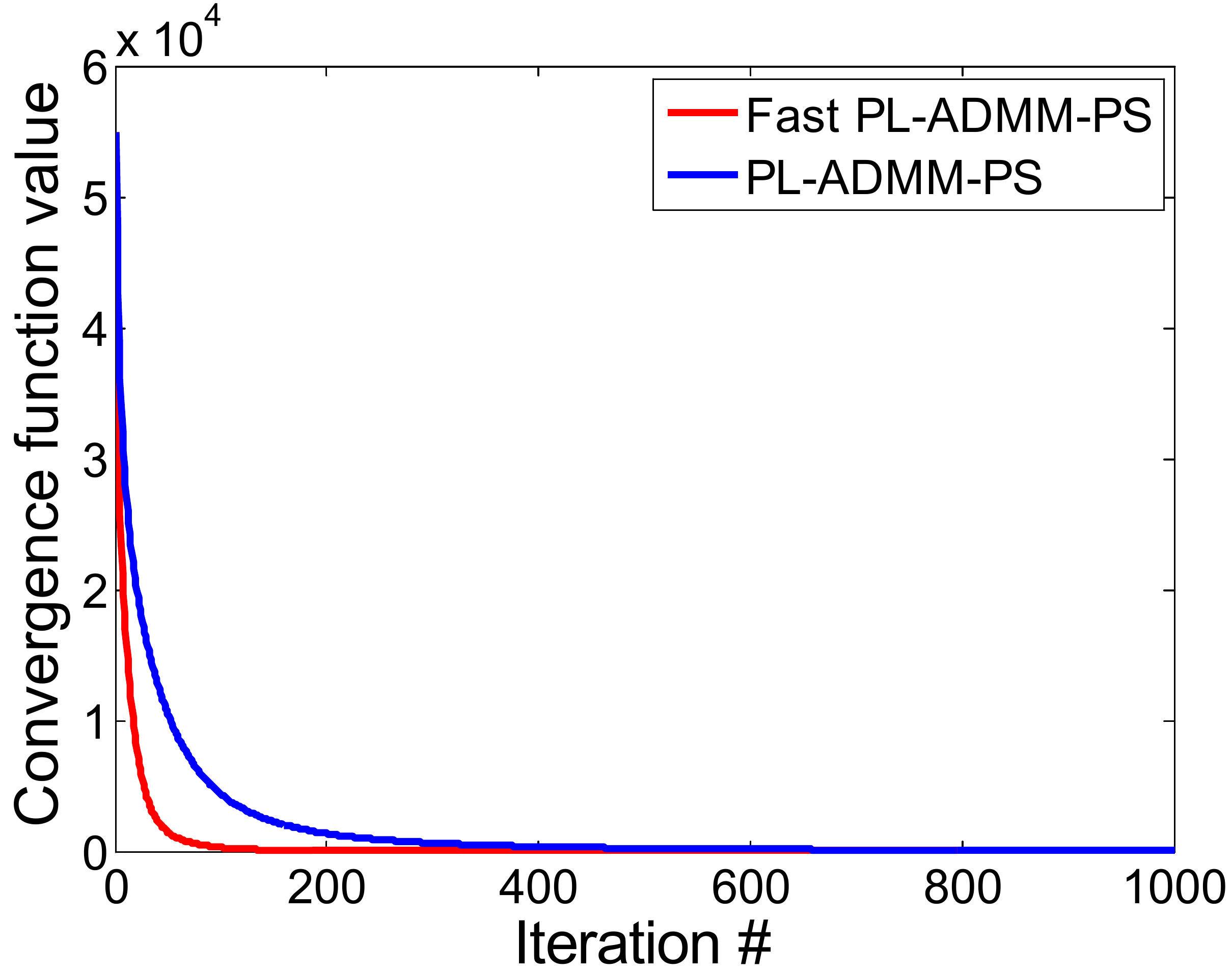}}
	\subfigure[$m=500$]{
		\includegraphics[width=0.32\textwidth]{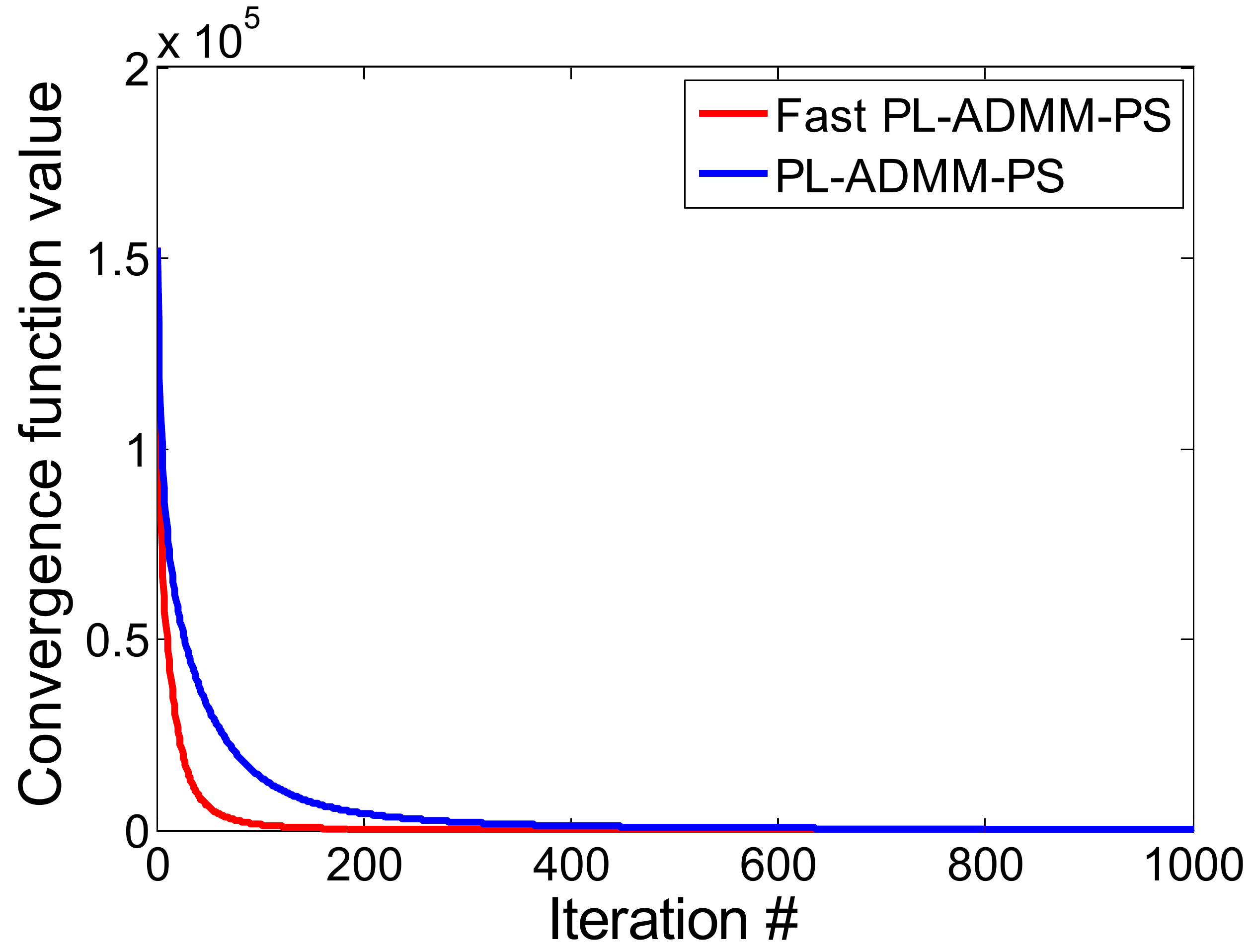}}
	\vspace{-0.2cm}
	\caption{Plots of the convergence function values of (\ref{con_FastPS}) in each iterations by using PL-ADMM-PS and Fast PL-ADMM-PS for (\ref{eqmulti}) with different sizes of $\X\in\mathbb{R}^{m\times m}$.} 
	\label{fig_res_fALMps}
\end{figure*}

\begin{figure*}[t]
	\setcounter{subfigure}{0}
	\subfigure[5 subjects]{
		\includegraphics[width=0.32\textwidth]{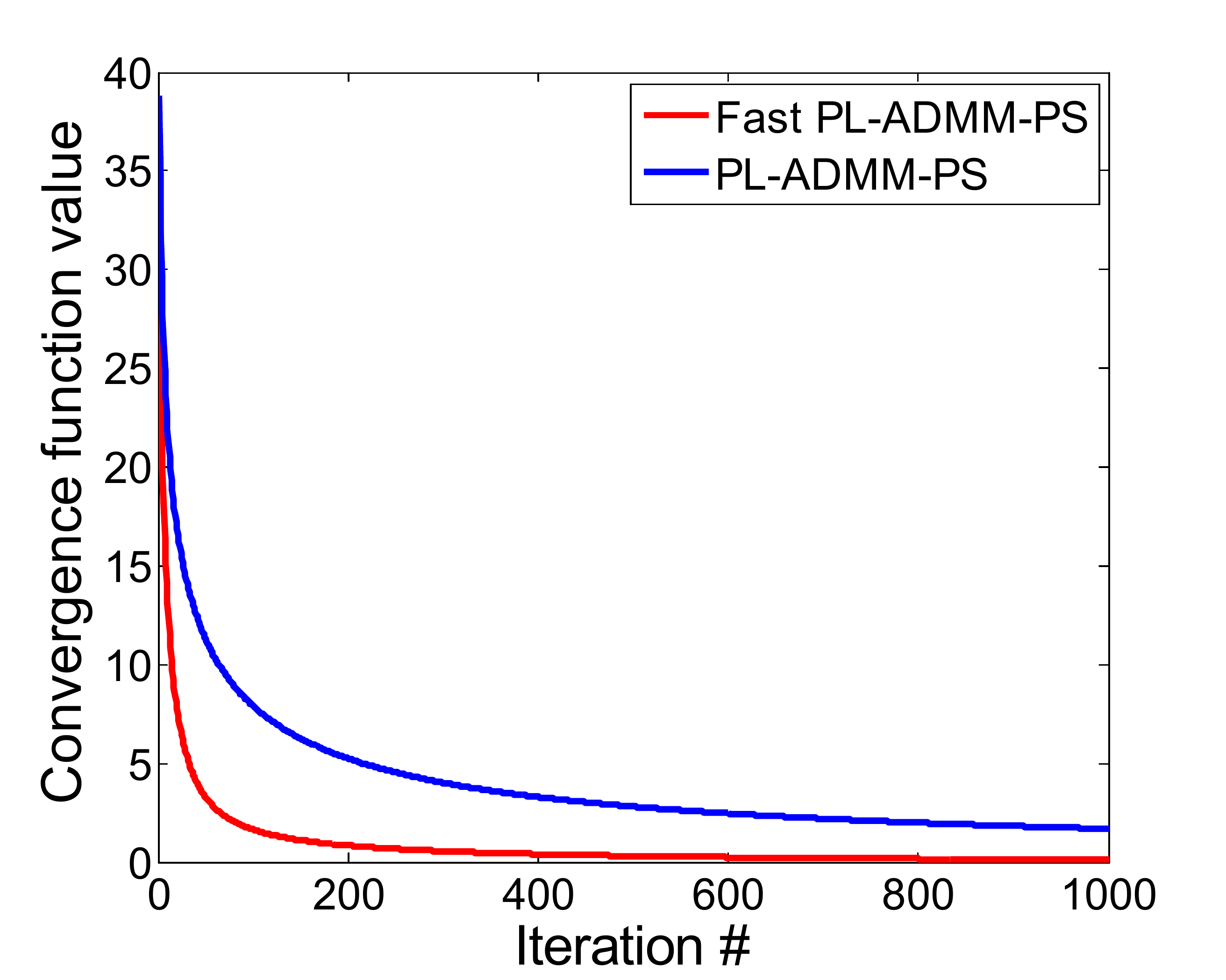}}
	\subfigure[8 subjects]{
		\includegraphics[width=0.32\textwidth]{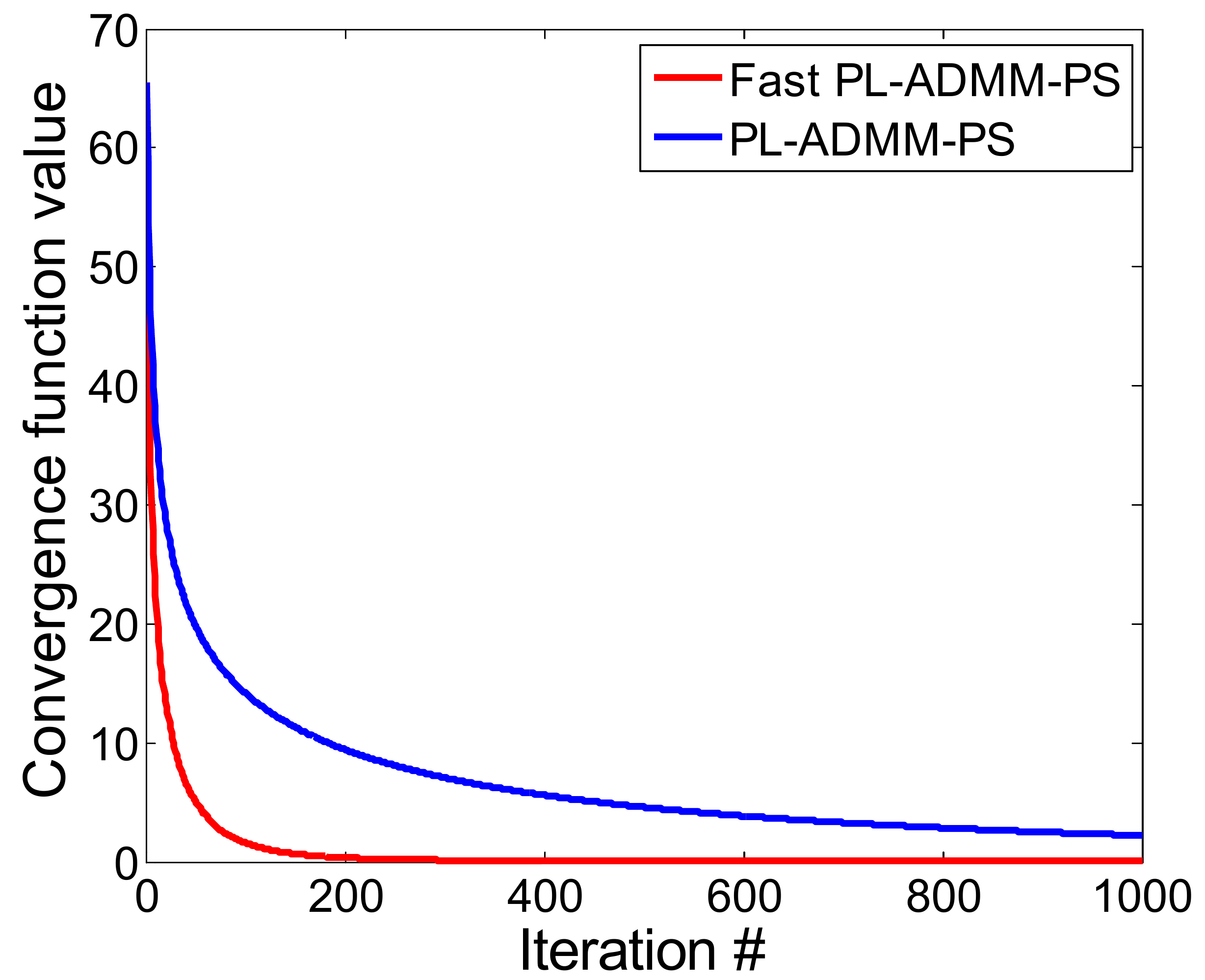}}
	\subfigure[10 subjects]{
		\includegraphics[width=0.32\textwidth]{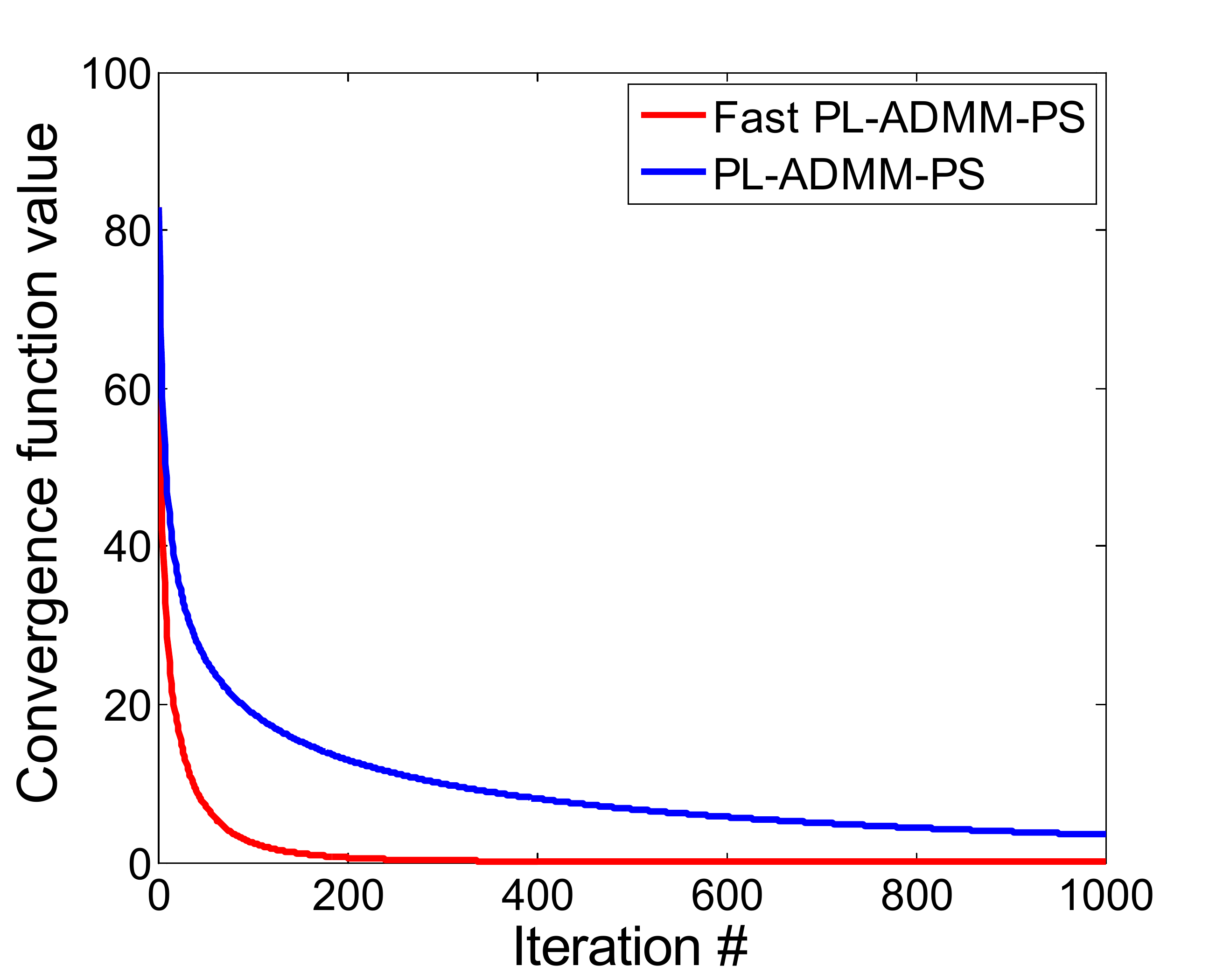}}
	\caption{Plots of the convergence function values of (\ref{con_FastPS}) in each iterations by using PL-ADMM-PS and Fast PL-ADMM-PS for (\ref{eqsub}) with different sizes of data  $\X$ for subspace clustering.} 
	\label{fig_res_ALMpsclu}
	\vspace{-0.5cm}
\end{figure*}

\subsection{Comparison of PL-ADMM-PS and Fast PL-ADMM-PS}
In this subsection, we conduct a problem with three blocks of variables as follows 
\begin{equation}\label{eqmulti}
\begin{split}
\min_{\X_1,\X_2,\X_3} & \ \  \sum_{i=1}^{3} \left(||\X_i||_{\ell_i}+\frac{\alpha_i }{2}||\mathbf{C}_i\X_i-\D_i||_F^2\right), \\
\text{s.t.} & \ \ \sum_{i=1}^{3}\mathbf{A}_i\X_i=\mathbf{B},
\end{split}
\end{equation}
where $||\cdot||_{\ell_1}=||\cdot||_1$ is the $\ell_1$-norm, $||\cdot||_{\ell_2}=||\cdot||_*$ is the nuclear norm, and $||\cdot||_{\ell_3}=||\cdot||_{2,1}$ is the $\ell_{2,1}$-norm defined as the sum of the $\ell_2$-norm of each column of a matrix. We simply consider all the matrices with the same size ${\mathbf{A}}_i, {\mathbf{C}}_i, {\D}_i, \mathbf{B}, {\X}_i\in\mathbb{R}^{m\times m}$. The matrices ${\mathbf{A}}_i, {\mathbf{C}}_i, {\D}_i, i=1,2,3$, and $\mathbf{B}$ are generated by the Matlab command \mcode{randn}. We set the parameters $\alpha_1=\alpha_2=\alpha_3=0.1$. Problem (\ref{eqmulti}) can be solved by PL-ADMM-PS and Fast PL-ADMM-PS, which have the same and cheap per-iteration cost. The experiments are conducted on three different values of $m=$100, 300 and 500. Figure \ref{fig_res_fALMps} plots the convergence function values of PL-ADMM-PS  and Fast PL-ADMM-PS in (\ref{con_FastPS}). It can be seen that Fast PL-ADMM-PS converges much faster than PL-ADMM-PS. Though Fast PL-ADMM-PS only accelerates PL-ADMM-PS for the smooth parts $g_i$'s, the improvement of Fast PL-ADMM-PS over PL-ADMM-PS is similar to that in Fast PALM over PALM. The reason behind this is that the Lipschitz constants $L_i$'s are not very small (around 400, 1200, and 2000 for the cases $m=100$, $m=300$, and $m=500$, respectively). And thus reducing the first term of (\ref{con_FastPS}) faster by our method is important.

\subsection{Application to Subspace Clustering}
In this subsection, we consider the following low rank and sparse representation problem for subspace clustering 
\begin{equation}\label{eqsub}
\begin{split}
\min_{\Z}  & \ \  \alpha_1||\Z||_*+\alpha_2||\Z||_1+\frac{1}{2}||\X\Z-\X||^2,\\
\text{s.t.}&\ \  \mathbf{1}^T\Z=\mathbf{1}^T,
\end{split}
\end{equation}
where $\X$ is the given data matrix. The above model is motivated by \cite{zhuang2012non}. However, we further consider the affine constraint $\mathbf{1}^T\Z=\mathbf{1}^T$ for affine subspace clustering \cite{6482137}. Problem (\ref{eqsub}) can be reformulated as a special case of problem (\ref{eq:model_problem_multivar}) by introducing auxiliary variables. Then it can be solved by PL-ADMM-PS and Fast PL-ADMM-PS. 

Given a data matrix $\X$ with each column as a sample, we solve (\ref{eqsub}) to get the optimal solution $\Z^*$. Then the affinity matrix $\W$ is defined as $\W=(|\Z|+|\Z^T|)/2$. The normalized cut algorithm \cite{shi2000normalized} is then performed on $\W$ to get the clustering results of the data matrix $\X$. The whole clustering algorithm is the same as \cite{6482137}, but using our defined affinity matrix $\W$ above.

\renewcommand{\arraystretch}{1.3}
\begin{table}[t] 
	\small
	\vspace{-0.2cm}
	\caption{Comparision of subspace clustering accuracies ($\%$) on the Extended Yale B database.}\vspace{0.1cm}
	\centering
	\begin{tabular}{c|c|c|c}
		\hline
		Methods  & 5 subjects & 8 subjects & 10 subjects \\ \hline
		PL-ADMM-PS & 94.06      & 85.94      & 75.31       \\ \hline
		Fast PL-ADMM-PS & 96.88      & 90.82      & 75.47       \\ \hline
	\end{tabular}\label{tabacc}
	\vspace{-0.5cm}
\end{table}

We conduct experiments on the Extended Yale B database \cite{georghiades2001few}, which is challenging for clustering. It consists of 2,414 frontal face images of 38 subjects under various lighting, poses and illumination conditions. Each subject has 64 faces.  We construct
three matrices $\X$ based on the first 5, 8 and
10 subjects. The data matrices $\X$ are first
projected into a $5\times 6$, $8\times 6$, and $10\times 6$-dimensional subspace
by PCA, respectively. Then we run PL-ADMM-PS and Fast PL-ADMM-PS for 1000 iterations, and use the solutions $\Z$ to define the affinity matrix $\W=(|\Z|+|\Z^T|)/2$. Finally, we can obtain the clustering results by normalized cuts.  The accuracy, calculated by the best matching rate of the predicted label and the ground
truth of data, is reported to measure the performance. Table \ref{tabacc} shows the clustering accuracies based on the solutions to problem (\ref{eqsub}) obtained by PL-ADMM-PS and Fast PL-ADMM-PS. It can be seen that Fast PL-ADMM-PS usually outperfoms PL-ADMM-PS since it achieves a better solution than PL-ADMM-PS within 1000 iterations. This can be verified in Figure \ref{fig_res_ALMpsclu} which shows the convergence function values  in (\ref{con_FastPS}) of PL-ADMM-PS and Fast PL-ADMM-PS in each iteration. It can be seen that our Fast PL-ADMM-PS converges much faster than PL-ADMM-PS.

\section{Conclusions}
This paper presented two fast solvers for the linearly constrained convex problem (\ref{eq:model_problem_multivar}). In particular, we proposed the Fast Proximal Augmented Lagragian Method (Fast PALM) which achieves the  convergence rate $O(1/K^2)$. Note that such a rate is theoretically optimal by comparing with the rate $O(1/K)$ by traditional ALM/PALM. Our fast version does not require additional assumptions (e.g. boundedness of $X$ and $\Lambda$, or a predefined number of iterations) as in the previous works~\cite{azadi2014towards,ouyang2014accelerated}. In order to further reduce the per-iteration complexity and handle the multi-blocks problems ($n>2$), we proposed the Fast Proximal Linearized ADMM with Parallel Splitting (Fast PL-ADMM-PS). It also achieves the optimal $O(1/K^2)$ rate for the smooth part of the objective. Compared with PL-ADMM-PS, though Fast PL-ADMM-PS requires additional assumptions on the boundedness of $X$ and $\Lambda$ in theory, our experimental results show that significant improvements are obtained especially when the Lipschitz constant of the smooth part is relatively large.

\section*{Acknowledgements}
This research is supported by the Singapore National Research Foundation under its International Research Centre
@Singapore Funding Initiative and administered by the IDM Programme Office. Z. Lin is supported by NSF China
(grant nos. 61272341 and 61231002), 973 Program of China
(grant no. 2015CB3525) and MSRA Collaborative Research
Program. 

{
	\footnotesize
	\bibliographystyle{aaai}
	\bibliography{FPLALM}
}

\newpage
\onecolumn
~\\\\
\begin{center}
	\LARGE\textbf{Supplementary Material of Fast Proximal Linearized Alternating Direction Method of Multiplier with Parallel Splitting}
\end{center}
\vspace{1em}
\renewcommand{\thefootnote}{$\star$} 
\begin{center}
	{\Large\textbf{Canyi Lu}$^1$, \textbf{Huan Li}$^2$, \textbf{Zhouchen Lin}$^{2,3,}$\footnote{Corresponding author.}, \textbf{Shuicheng Yan}$^1$}\\
	$^1$ Department of Electrical and Computer Engineering, National University of Singapore\\
	$^2$ Key Laboratory of Machine Perception (MOE), School of EECS, Peking University\\
	$^3$ Cooperative Medianet Innovation Center, Shanghai Jiaotong University \\
	{\tt\small canyilu@gmail.com,  lihuanss@pku.edu.cn, zlin@pku.edu.cn, eleyans@nus.edu.sg}
\end{center}

This documents provides the proof details of the convergence results of our proposed fast methods. First, in Section 1, we give some useful results which are useful for the convergence analysis of Fast PALM in Section 2 and Fast PL-ADMM-PS in Section 3.

\section{Some Lemmas}
\label{sec1}
\begin{lemma}\label{Lem_lips}
	Let $g: \mathbb{R}^{m}\rightarrow\mathbb{R}$ be a continuously differentiable function with Lipschitz continuous gradient and Lipschitz constant $L$. Then, for any ${\x}, {\y}\in\mathbb{R}^m$,
	\begin{equation}\label{lippro}
	g({\x})\leq g({\y})+\langle {\x}-{\y},\nabla g({\y}) \rangle+\frac{L}{2}||{\x}-{\y}||^2.
	\end{equation}
\end{lemma}
\begin{lemma}\label{Lem_Pythagoras}
	Given any $\bm{a}$, $\bm{b}$, $\bm{c}$, $\bm{d} \in \mathbb{R}^m $, we have
	\begin{equation}\label{lempytha1}
	\langle\bm{a}-\bm{b},\bm{a}-\bm{c}\rangle=\frac{1}{2}\left(\| \bm{a}-\bm{b} \|^2+\|\bm{a}-\bm{c}\|^2-\|\bm{b}-\bm{c}\|^2\right).
	\end{equation}
	\begin{equation}\label{lempytha2}
	\langle\bm{a}-\bm{b},\bm{c}-\bm{d}\rangle=\frac{1}{2}\left(\| \bm{a}-\bm{d} \|^2-\|\bm{a}-\bm{c}\|^2-\|\bm{b}-\bm{d}\|^2+\|\bm{b}-\bm{c}\|^2\right).
	\end{equation}
\end{lemma}
\begin{lemma}
	Assume the sequences $\{\ak\}$ and $\{\bk\}$ satisfy $a^{(0)}=1$, $0<\akk-\ak\leq1$ and $\bk>0$. Then we have
	\begin{equation}\label{lemseq333}
	\sumk \ak(\bk-\bkk)\leq \sumk\bk.
	\end{equation}
\end{lemma}
\emph{Proof.} We deduce
\begin{eqnarray*}
	\sumk \ak(\bk-\bkk)&=& a^{(0)}b^{(0)}+\sum_{k=0}^{K-1}(\akk-\ak)\bkk-a^{(K)}b^{(K+1)} \\
	&\leq&b^{(0)}+\sum_{k=0}^{K-1}\bkk=\sumk\bk.
\end{eqnarray*}
$\hfill\blacksquare$
\begin{lemma}\label{Lem_theta}
	Define the sequence $\{\thetak\}$ as $\theta^{(0)}=1$, $\frac{1-\thetakk}{(\thetakk)^2}=\frac{1}{(\thetak)^2}$ and $\thetak>0$. Then we have the following properties
	\begin{equation}
	\thetakk=\frac{-(\thetak)^2+\sqrt{(\thetak)^4+4(\thetak)^2}}{2},
	\end{equation}
	\begin{equation}\label{lem3pro1}
	\sum\limits_{k=0}^K \frac{1}{\thetak}=\frac{1}{(\theta^{(K)})^2},
	\end{equation}
	\begin{equation}\label{lem3pro2}
	0<\frac{1}{\thetakk}-\frac{1}{\thetak}<1,
	\end{equation}
	\begin{equation}\label{lem3pro3}
	\thetak\leq\frac{2}{k+2},
	\end{equation}
	and
	\begin{equation}\label{lem3pro4}
	\thetak\leq1.
	\end{equation}
\end{lemma}
\emph{Proof.} From the definition of $\thetakk$, it is easy to get that $\thetakk=\frac{-(\thetak)^2+\sqrt{(\thetak)^4+4(\thetak)^2}}{2}$. This implies that $\thetak$ is well defined for any $k\geq0$. Furthermore, since $\frac{1}{\thetakk}=\frac{1}{(\thetakk)^2}-\frac{1}{(\thetak)^2}$ and $\theta^{(0)}=1$, we have
\begin{eqnarray}
\sum\limits_{k=0}^K \frac{1}{\thetak}=\frac{1}{\theta^{(0)}}+\sum\limits_{k=0}^{K-1} \frac{1}{\thetakk}=\frac{1}{\theta^{(0)}}+\sum\limits_{k=0}^{K-1} \left(\frac{1}{(\thetakk)^2}-\frac{1}{(\thetak)^2}\right)=\frac{1}{\theta^{(0)}}+\frac{1}{(\theta^{(K)})^2}-\frac{1}{(\theta^{(0)})^2}    =\frac{1}{(\theta^{(K)})^2}.
\end{eqnarray}
From $\frac{1}{\thetakk}=\frac{1}{(\thetakk)^2}-\frac{1}{(\thetak)^2}$, $\thetak>0$ and $\theta^{(k-1)}>0$, we can easily get
\begin{equation}
\frac{1}{\thetakk}-\frac{1}{\thetak}> 0,
\end{equation}
and
\begin{equation}
\frac{1}{\thetakk}-\frac{1}{\thetak}=\frac{1}{\thetakk}-\frac{\sqrt{1-\thetakk}}{\thetakk}=\frac{1-\sqrt{1-\thetakk}}{\thetakk}=\frac{1}{1+\sqrt{1-\thetakk}}<1.
\end{equation}
Next we proof $\thetak\leq\frac{2}{k+2}$ by induction. First $\theta^{(0)}=1\leq\frac{2}{0+2}$. Now assume that $\thetak\leq\frac{2}{k+2}$ and we prove $\thetakk\leq\frac{2}{k+3}$. We deduce
\begin{equation*}
\begin{split}
\thetakk&=\frac{-(\thetak)^2+\sqrt{(\thetak)^4+4(\thetak)^2}}{2}=\frac{2(\thetak)^2}{(\thetak)^2+\sqrt{(\thetak)^4+4(\thetak)^2}}\\
&=\frac{2}{1+\sqrt{1+\frac{4}{(\thetak)^2}}}\leq \frac{2}{1+\sqrt{1+(k+2)^2}}\leq \frac{2}{k+3}.
\end{split}
\end{equation*}
So $(\ref{lem3pro3})$ holds. Note that $\thetak$ is decreasing by (\ref{lem3pro2}) and  $\theta^{(0)}=1$, we have (\ref{lem3pro4}). The proof is completed.
$\hfill\blacksquare$

\section{Convergence Analysis of Fast PALM}
\label{sec2}
In this section, we give the convergence analysis of Fast PALM for solving the following problem
\begin{equation}
\min\limits_{\x} f(\x),\quad
s.t.\quad
\A(\x)=\b,\label{eq:model_problem_onevar}
\end{equation}
where $f(\x)=g(\x)+h(\x)$, both $g$ and $h$ are convex,
and $g \in C^{1.1}$:
\begin{equation}
\lbar \nabla g(\x)- \nabla g(\y)\rbar \leq L \lbar \x-\y\rbar,\quad \forall \x,\y. \label{eq:Lipschitz-continue}
\end{equation}
For the completeness, we give the Fast PALM in Algorithm \ref{alg_fastpalm}. 

It is worth mentioning that the definition of $\thetakk$ in (\ref{defthetakk1}) is equivalent to $\theta^{(0)}=1$, $\frac{1-\thetakk}{(\thetakk)^2}=\frac{1}{(\thetak)^2}$ and $\thetak>0$ in Lemma \ref{Lem_theta}. Such a property will be used in the following analysis several times.

The anaysis of our algorithms is based on the following property:
\begin{lemma}\label{Lemma_opt}
	$\tilde{\x}$ is an optimal solution to
	(\ref{eq:model_problem_onevar}) if and only if there exists
	$\alpha > 0$, such that
	\begin{equation}
	f(\tilde{\x})-f(\x^*)+\langle\bflambda^*,\A(\tilde{\x})-\b\rangle+\frac{\alpha}{2}\lbar\A(\tilde{\x})-\bm{b}\rbar^2=0.
	\end{equation}
\end{lemma}

\begin{algorithm}[!htp]\label{alg_fastpalm}
	\caption{Fast PALM Algorithm}
	\hrule
	\vspace{0.1cm}
	\textbf{Initialize}: $\x^0$, $\z^0$, $\bm{\lambda}^0$, $\beta^{(0)}=\theta^{(0)}=1$.\\
	\For{$k = 0, 1, 2,\cdots$}{
		\begin{eqnarray}
		\ykk&=&(1-\thetak)\x^k+\thetak\zk;\\
		\z^{k+1}&=&\argmin\limits_{\x} g(\y^{k+1})+ \left \langle\nabla g(\y^{k+1}),\x-\y^{k+1} \right \rangle+h(\x)\notag\\
		&&+ \left \langle\bm\lambda^k,\A(\x)-\b \right \rangle+\frac{\betak}{2}\|\A(\x)-\b\|^2+\frac{L\thetak}{2}\|\x-\z^k\|^2;\label{z_update_onevar}\\
		\x^{k+1}&=&(1-\thetak)\x^k+\thetak\z^{k+1}; \label{ddefinex}\\
		\bm{\lambda}^{k+1}&=&\bm{\lambda}^k+\betak(\A(\z^{k+1})-\b);\label{update_lambda_one} \\
		\thetakk&=&\frac{-(\thetak)^2+\sqrt{(\thetak)^4+4(\thetak)^2}}{2}; \label{defthetakk1}\\
		\betakk&=& \frac{1}{\thetakk}.
		\end{eqnarray}
	}
	\hrule
	\vspace{0.1cm}
\end{algorithm}

\begin{proposition}\label{prop_one}
	In Algorithm \ref{alg_fastpalm},  for any $\x$, we have
	\begin{eqnarray}
	&&\frac{1-\thetakk}{(\thetakk)^2} \left(f(\x^{k+1})-f(\x)\right)-\frac{1}{\thetak}\left \langle\A^T(\blambdakk),\x-\z^{k+1} \right \rangle\qquad \notag\\
	&\leq&\frac{1-\thetak}{(\thetak)^2}\left(f(\x^k)-f(\x)\right) + \frac{L}{2}\left(\|\z^{k}-\x\|^2-\|\z^{k+1}-\x\|^2\right).\label{eqpro1}
	\end{eqnarray}
\end{proposition}
\emph{Proof.}
From the optimality of $\zkk$ to (\ref{z_update_onevar}), we have
\begin{eqnarray}
0&\in& \partial h(\z^{k+1})+\nabla g(\y^{k+1})+\A^T(\bm{\lambda}^k)+\betak\A^T(\A(\z^{k+1})-\b)+L\thetak(\z^{k+1}-\z^k)\notag\\
&=&\partial h(\z^{k+1})+\nabla g(\y^{k+1})+\A^T(\bm{\lambda}^{k+1})+L\thetak(\z^{k+1}-\z^k),\label{pro1proof199}
\end{eqnarray}
where (\ref{pro1proof199}) uses (\ref{update_lambda_one}).
From the convexity of $h$, we have
\begin{equation}
h(\x)-h(\z^{k+1})\geq  \left \langle-\nabla g(\y^{k+1})-\A^T(\bm{\lambda}^{k+1})-L\thetak(\z^{k+1}-\z^k),\x-\z^{k+1} \right \rangle. \label{eq:opt_z2_one}
\end{equation}
On the other hand,
\begin{eqnarray}
f(\x^{k+1})&\leq& g(\y^{k+1})+ \left \langle\nabla g(\y^{k+1}),\x^{k+1}-\y^{k+1} \right \rangle+\frac{L}{2}\|\x^{k+1}-\y^{k+1}\|^2+h(\x^{k+1})\label{t2l1}\\
&=&g(\y^{k+1})+ \left \langle\nabla g(\y^{k+1}),(1-\thetak)\x^k+\thetak\z^{k+1}-\y^{k+1} \right \rangle\notag\\
&&+\frac{L}{2}\|(1-\thetak)\x^k+\thetak\z^{k+1}-\y^{k+1}\|^2+h\left((1-\thetak)\x^k+\thetak\z^{k+1}\right)\label{t2l2}\\
&\leq&(1-\thetak)\left(g(\y^{k+1})+ \left \langle\nabla g(\y^{k+1}),\x^k-\y^{k+1} \right \rangle+h(\x^k)\right)\notag\\
&&+\thetak\left(g(\y^{k+1})+ \left \langle\nabla g(\y^{k+1}),\z^{k+1}-\y^{k+1} \right \rangle+h(\z^{k+1})\right)+\frac{L(\thetak)^2}{2}\|\z^{k+1}-\z^k\|^2\label{t2l3}\\
&=&(1-\thetak)\left(g(\y^{k+1})+ \left \langle\nabla g(\y^{k+1}),\x^k-\y^{k+1} \right \rangle+h(\x^k)\right)\notag\\
&&+\thetak\left(g(\y^{k+1})+ \left \langle\nabla g(\y^{k+1}),\x-\y^{k+1} \right \rangle+ \left \langle\nabla g(\y^{k+1}),\z^{k+1}-\x \right \rangle+h(\z^{k+1})\right)\notag\\
&&+\frac{L(\thetak)^2}{2}\|\z^{k+1}-\z^k\|^2\notag\\
&\leq&(1-\thetak)f(\x^k)+\thetak\left(g(\x)+ \left \langle\nabla g(\y^{k+1}),\z^{k+1}-\x \right \rangle+h(\z^{k+1})\right)+\frac{L(\thetak)^2}{2}\|\z^{k+1}-\z^k\|^2\label{t2l5}\\
&\leq&(1-\thetak)f(\x^k)+\thetak\left(g(\x)+h(\x)+ \left \langle\A^T(\bm{\lambda}^{k+1})+L\thetak(\z^{k+1}-\z^k),\x-\z^{k+1} \right \rangle\right)\label{t2l6} \\
&&+\frac{L(\thetak)^2}{2}\|\z^{k+1}-\z^k\|^2\notag\\
&=&(1-\thetak)f(\x^k)+\thetak f(\x)+\thetak \left \langle\A^T(\bm{\lambda}^{k+1}),\x-\z^{k+1} \right \rangle
-\frac{L(\thetak)^2}{2}\left(\|\z^{k+1}-\x\|^2-\|\z^k-\x\|^2\right),\label{t2l7}
\end{eqnarray}
where (\ref{t2l1}) uses (\ref{lippro}), (\ref{t2l2}) uses (\ref{ddefinex}),  (\ref{t2l3}) is from the convexity of $h$, (\ref{t2l5}) is from the convexity of $g$, (\ref{t2l6}) uses (\ref{eq:opt_z2_one}) and (\ref{t2l7}) uses (\ref{lempytha1}).
Rerangging the above inequality leads to
\begin{eqnarray}
&&\left(f(\x^{k+1})-f(\x)\right)-\thetak\left \langle\A^T(\blambdakk),\x-\z^{k+1} \right \rangle\qquad \notag\\
&\leq&(1-\thetak)\left(f(\x^k)-f(\x)\right) + \frac{L(\thetak)^2}{2}\left(\|\z^{k}-\x\|^2-\|\z^{k+1}-\x\|^2\right).\label{eq299}
\end{eqnarray}
Diving both sides of the above inequality by $(\thetak)^2$ leads to
\begin{eqnarray*}
	&&\frac{1}{(\thetak)^2} \left(f(\x^{k+1})-f(\x)\right)-\frac{1}{\thetak}\left \langle\A^T(\blambdakk),\x-\z^{k+1} \right \rangle\qquad \notag\\
	&\leq&\frac{1-\thetak}{(\thetak)^2}\left(f(\x^k)-f(\x)\right) + \frac{L}{2}\left(\|\z^{k}-\x\|^2-\|\z^{k+1}-\x\|^2\right).
\end{eqnarray*}
The proof is completed by using the property of $\thetak$ in Lemma \ref{Lem_theta}. $\hfill\blacksquare$

\begin{proposition}\label{prop2}
	In Algorithm \ref{alg_fastpalm}, the following result holds for any $\blambda$
	\begin{eqnarray}
	&& \langle\A(\zkk)-\b,\blambda-\blambdakk\rangle+\frac{\betak}{2}\|\A(\zkk)-\b\|^2 \notag\\
	& = &  \frac{1}{2\betak} \left(\|\blambdak-\blambda\|^2-\|\blambdakk-\blambda\|^2\right)\label{eqpro2}
	\end{eqnarray}
\end{proposition}
\emph{Proof.} By using (\ref{update_lambda_one}) and (\ref{lempytha1}), we have
\begin{eqnarray*}
	&& \langle\A(\zkk)-\b,\blambda-\blambdakk\rangle \\
	&=& \frac{1}{\betak}\langle\blambdakk-\blambdak,\blambda-\blambdakk\rangle \\
	&=& \frac{1}{2\betak} \left(\|\blambdak-\blambda\|^2-\|\blambdakk-\blambda\|^2-\|\blambdakk-\blambdak\|^2 \right) \\
	&=& \frac{1}{2\betak} \left(\|\blambdak-\blambda\|^2-\|\blambdakk-\blambda\|^2\right)- \frac{\betak}{2}\|\A(\zkk)-\b\|^2.
\end{eqnarray*}
The proof is completed. $\hfill\blacksquare$

\begin{theorem}\label{conv_rate_one}
	In Algorithm \ref{alg_fastpalm}, for any $K>0$, we have
	\begin{eqnarray}
	&&f(\x^{K+1})-f(\x^*)+\left \langle\bm\lambda^*,\A(\x^{K+1})-\b \right \rangle+\frac{1}{2}\|\A (\x^{K+1})-\b\|^2\qquad\\
	&\leq&\frac{2}{(K+2)^2}\left(LD_{\x^*}^2+D_{\lambda^*}^2\right).
	\end{eqnarray}
\end{theorem}
\emph{Proof.} Let $\x=\x^*$ and $\blambda=\blambda^*$ in (\ref{eqpro1}) and (\ref{eqpro2}). We have
\begin{eqnarray}
&&\frac{1-\thetakk}{(\thetakk)^2}\left(f(\x^{k+1})-f(\x^*)\right)-\frac{1-\thetak}{(\thetak)^2}\left(f(\x^k)-f(\x^*)\right)+\frac{1}{\thetak}\left \langle\bm\lambda^*,\A(\z^{k+1})-\b \right \rangle\qquad\label{eqproof399}\\
&\leq&\frac{L}{2}\left(\|\z^{k}-\x^*\|^2-\|\z^{k+1}-\x^*\|^2\right)+\frac{1}{\thetak}\langle\blambda^*-\blambdakk,\A(\zkk)-\b\rangle \label{eqproof40} \\
&\leq&\frac{L}{2}\left(\|\z^{k}-\x^*\|^2-\|\z^{k+1}-\x^*\|^2\right)+\frac{1}{2\thetak\betak}\left( \|\bm{\lambda}^k-\bm\lambda^*\|^2-\|\bm{\lambda}^{k+1}-\bm\lambda^*\|^2 \right)-\frac{\betak}{2\thetak}\|\A( \z^{k+1})-\b\|^2\label{eqproof41} \\
&=& \frac{L}{2}\left(\|\z^{k}-\x^*\|^2-\|\z^{k+1}-\x^*\|^2\right)+\frac{1}{2}\left( \|\bm{\lambda}^k-\bm\lambda^*\|^2-\|\bm{\lambda}^{k+1}-\bm\lambda^*\|^2 \right)-\frac{1}{2(\thetak)^2}\|\A( \z^{k+1})-\b\|^2,\label{eqproof42}
\end{eqnarray}
where (\ref{eqproof40}) uses the fact $\A(\x^*)=\b$, (\ref{eqproof41}) uses (\ref{eqpro2}) and (\ref{eqproof42}) uses $\betak=\frac{1}{\thetak}$.

Summing (\ref{eqproof399})-(\ref{eqproof42}) from $k=0$ to $K$, we have
\begin{eqnarray}
&&\frac{1-\theta^{(K+1)}}{(\theta^{(K+1)})^2}\left(f(\x^{K+1})-f(\x^*)\right)-\frac{1-\theta^{(0)}}{(\theta^{(0)})^2}\left(f(\x^{0})-f(\x^*)\right)+\sum_{k=0}^K\frac{1}{\thetak}\left \langle\bm\lambda^*,\A(\z^{k+1})-\b \right \rangle\label{pro1eq4441}\\
&\leq&\frac{L}{2}\|\z^{0}-\x^*\|^2+\frac{1}{2}\|\bm\lambda^{0}-\bm\lambda^*\|^2-\sum_{k=0}^K\frac{1}{2(\thetak)^2}\|\A (\z^{k+1})-\b\|^2\notag\\
&\leq&\frac{L}{2}\|\z^{0}-\x^*\|^2+\frac{1}{2}\|\bm\lambda^{0}-\bm\lambda^*\|^2-\sum_{k=0}^K\frac{1}{2\thetak}\|\A (\z^{k+1})-\b\|^2,\label{pro1eq444}
\end{eqnarray}
where (\ref{pro1eq444}) uses (\ref{lem3pro4}). Also note that $\theta^{(0)}=1$. So the second term of (\ref{pro1eq4441}) disappears.

On the other hand, by the property of $\thetak$ in Lemma \ref{Lem_theta} and (\ref{ddefinex}), we have
\begin{eqnarray}
&& \sum_{k=0}^K \frac{\zkk}{\thetak} \notag \\
&=&\sum_{k=0}^K\left(\frac{1}{(\thetak)^2}\xkk-\frac{1-\thetak}{(\thetak)^2}\xk\right)\notag  \\
&=&\sum_{k=0}^K\left(\frac{1-\thetakk}{(\thetakk)^2}\xkk-\frac{1-\thetak}{(\thetak)^2}\xk\right) \notag \\
&=&\frac{1-\theta^{(K+1)}}{(\theta^{(K+1)})^2}\x^{K+1}-\frac{1-\theta^{(0)}}{(\theta^{(0)})^2}\x^0 \notag \\
&=&\frac{1-\theta^{(K+1)}}{(\theta^{(K+1)})^2}\x^{K+1}\notag  \\
&=&\frac{1}{(\theta^{(K)})^2}\x^{K+1}.\label{profpro1eq455666}
\end{eqnarray}
So
\begin{eqnarray}\label{pro1proofeq46666}
\sum_{k=0}^K\frac{1}{\thetak}\left \langle\bm\lambda^*,\A(\z^{k+1})-\b \right \rangle=\frac{1}{(\theta^{(K)})^2}\left \langle\bm\lambda^*,\A(\x^{K+1})-\b \right \rangle.
\end{eqnarray}
By the convexity of $\|\cdot\|^2$, we have
\begin{eqnarray}
&&\sum_{k=0}^K\frac{1}{2\thetak}\|\A (\z^{k+1})-\b\|^2\notag\\
&=&\frac{1}{2(\theta^{(K)})^2}  \sum_{k=0}^K\frac{\theta^{(K)})^2}{\thetak}\|\A (\z^{k+1})-\b\|^2\notag\\
&\geq&\frac{1}{2(\theta^{(K)})^2}\left\|(\theta^{(K)})^2 \A\left(\sum_{k=0}^{K}\frac{\zkk}{\thetak}  \right)  -\b \right\|^2 \label{pro1proofeq488}\\
&=&\frac{1}{2(\theta^{(K)})^2}\|\A (\x^{K+1})-\b\|^2\label{pro1proofeq499},
\end{eqnarray}
where (\ref{pro1proofeq488}) uses (\ref{lem3pro1})  and (\ref{pro1proofeq499}) uses  (\ref{profpro1eq455666}).

Substituting (\ref{pro1proofeq46666}) into (\ref{pro1eq4441}) and (\ref{pro1proofeq499}) into (\ref{pro1eq444}) respectively, we obtain
\begin{eqnarray}
&&\frac{1-\theta^{(K+1)}}{(\theta^{(K+1)})^2}\left(f(\x^{K+1})-f(\x^*)\right)+\frac{1}{(\theta^{(K)})^2}\left \langle\bm\lambda^*,\A(\x^{K+1})-\b \right \rangle+\frac{1}{2(\theta^{(K)})^2}\|\A (\x^{K+1})-\b\|^2\qquad\label{pro1proofeq5001}\\
&\leq&\frac{L}{2}\|\z^{0}-\x^*\|^2+\frac{1}{2}\|\bm\lambda^{0}-\bm\lambda^*\|^2.\label{pro1proofeq51111}
\end{eqnarray}
Multiplying (\ref{pro1proofeq5001}) and (\ref{pro1proofeq51111}) by $(\theta^{(K)})^2$ and using (\ref{lem3pro3}) leads to
\begin{eqnarray*}
	&&f(\x^{K+1})-f(\x^*)+\left \langle\bm\lambda^*,\A(\x^{K+1})-\b \right \rangle+\frac{1}{2}\|\A (\x^{K+1})-\b\|^2\qquad\\
	&\leq&\frac{2}{(K+2)^2}\left(L\|\z^{0}-\x^*\|^2+\|\bm\lambda^{0}-\bm\lambda^*\|^2\right)\\
	&=&\frac{2}{(K+2)^2}\left(LD_{\x^*}^2+D_{\lambda^*}^2\right).
\end{eqnarray*}
The proof is completed.
$\hfill\blacksquare$


\section{Convergence Analysis of Fast PL-ADMM-PS}
\label{sec3}
In this section, we give the convergence analysis of Fast PL-ADMM-PS for solving the following problem
\begin{equation}
\min\limits_{\x_1,\cdots,\x_n} \sum\limits_{i=1}^n f_i(\x_i),\quad
s.t.\quad
\sum\limits_{i=1}^n\A_i(\x_i)=\b,\label{eq:model_problem_multivar}
\end{equation}
where $f_i(\x_i)=g_i(\x_i)+h_i(\x_i)$, both $g_i$ and $h_i$ are convex,
and $g_i \in C^{1.1}$. The whole procedure of Fast PL-ADMM-PS is shown in Algorithm \ref{alg_fastplalmps}.
\begin{algorithm}[!t]\label{alg_fastplalmps}
	\DontPrintSemicolon
	\caption{Fast PL-ADMM-PS Algorithm}
	\hrule
	\hrule
	\hrule
	\vspace{0.1cm}
	\textbf{Initialize}: $\x^0$, $\z^0$, $\bm{\lambda}^0$, $\theta^{(0)}=1$, fix $\betak=\beta$ for  $k\geq0$, $\eta_i>n\|\A_i\|^2$, $i=1,\cdots,n$,\\
	\For{$k=0, 1, 2, \cdots$}{
		\hspace*{1.8cm}$\parallelsum$  Update $\y_i$, $\z_i$, $\x_i$, $i=1,\cdots,n$, in parallel by
		\begin{align}
		\y_i^{k+1}&=(1-\thetak)\x_i^k+\thetak\z_i^k;\label{updateyykkk}\\
		\z_i^{k+1}&=\argmin\limits_{\x_i}  \left \langle\nabla g_i(\y_i^{k+1}),\x_i \right \rangle +h_i(\x_i)+ \left \langle\bm{\lambda}^k,\A_i(\x_i) \right \rangle+\left\langle \betak\A_i^T\left( \A(\z^k)-\b\right),\x_i\right\rangle   \notag \\
		&+\frac{L(g_i)\thetak+\betak\eta_i}{2}\|\x_i-\z_i^k\|^2; \label{update_z}\\
		\x_i^{k+1}&=(1-\thetak)\x_i^k+\thetak\z_i^{k+1};  \label{updatexi}\\
		\mathbf{\bm{\lambda}}^{k+1}&=\mathbf{\bm{\lambda}}^k+\beta^k\left(\A(\z^{k+1})-\b\right);\label{update_lambda}\\
		\thetakk&=\frac{-(\thetak)^2+\sqrt{(\thetak)^4+4(\thetak)^2}}{2}.
		\end{align}
	}
	\hrule
	\hrule
	\hrule
	\vspace{0.1cm}
\end{algorithm}


\begin{proposition}\label{prop3}
	In Algorithm \ref{alg_fastplalmps}, for any $\x_i$, we have
	\begin{eqnarray}
	&&\frac{1-\thetakk}{(\thetakk)^2}\left(f_i(\x_i^{k+1})-f_i(\x_i)\right)-\frac{1}{\thetak} \left \langle\A_i^T({\hatlambdakk}),\x_i-\z_i^{k+1} \right \rangle\notag\\
	&\leq&\frac{1-\thetak}{(\thetak)^2}\left(f_i(\x_i^k)-f_i(\x_i)\right)+\frac{L_i}{2}\left(\|\z_i^{k}-\x_i\|^2-\|\z_i^{k+1}-\x_i\|^2\right)\notag\\
	&&+\frac{\betak\eta_i}{2\thetak}\left(\|\z_i^{k}-\x_i\|^2-\|\z_i^{k+1}-\x_i\|^2-\|\z_i^{k+1}-\z_i^k\|^2\right),\label{pro331}
	\end{eqnarray}
	where
	\begin{equation}\label{deflambdahat}
	\hat{\bm\bflambda}^{k+1}=\bm\lambda^k+\betak\left(\A(\z^{k})-\b\right).
	\end{equation}
\end{proposition}
\emph{Proof.} From the optimality of $\zkk_i$ to (\ref{update_z}), we have
\begin{eqnarray}
0&\in& \partial h_i(\z_i^{k+1})+\nabla g_i(\y_i^{k+1})+\A_i^T(\bm{\lambda}^k)+\betak\A_i^T(\A(\z^{k})-\b)+(L_i\thetak+\betak\eta_i)(\z_i^{k+1}-\z_i^k)\notag\\
&=&\partial h_i(\z_i^{k+1})+\nabla g_i(\y_i^{k+1})+\A_i^T(\bm{\hat\lambda}^{k+1})+(L_i\thetak+\betak\eta_i)(\z_i^{k+1}-\z_i^k),\label{pro353}
\end{eqnarray}
where (\ref{pro353}) uses (\ref{deflambdahat}). From the convexity of $h_i$, we have
\begin{equation}
h_i(\x_i)-h_i(\z_i^{k+1})\geq  \left \langle-\nabla g_i(\y_i^{k+1})-\A_i^T(\bm{\hat\lambda}^{k+1})-(L_i\thetak+\betak\eta_i)(\z_i^{k+1}-\z_i^k),\x_i-\z_i^{k+1} \right \rangle. \label{eq:opt_z2}
\end{equation}
On the other hand,
\begin{eqnarray}
f_i(\x_i^{k+1})&\leq& g_i(\y_i^{k+1})+ \left \langle\nabla g_i(\y_i^{k+1}),\x_i^{k+1}-\y_i^{k+1} \right \rangle+\frac{L_i}{2}\|\x_i^{k+1}-\y_i^{k+1}\|^2+h_i(\x_i^{k+1})\label{p1l1}\\
&=&g_i(\y_i^{k+1})+ \left \langle\nabla g_i(\y_i^{k+1}),(1-\thetak)\x_i^k+\thetak\z_i^{k+1}-\y_i^{k+1} \right \rangle\notag\\
&&+\frac{L_i}{2}\|(1-\thetak)\x_i^k+\thetak\z_i^{k+1}-\y_i^{k+1}\|^2+h_i\left((1-\thetak)\x_i^k+\thetak\z_i^{k+1}\right)\label{pro3566}\\
&\leq&(1-\thetak)\left(g_i(\y_i^{k+1})+ \left \langle\nabla g_i(\y_i^{k+1}),\x_i^k-\y_i^{k+1} \right \rangle+h_i(\x_i^k)\right)\notag\\
&&+\thetak\left(g_i(\y_i^{k+1})+ \left \langle\nabla g_i(\y_i^{k+1}),\z_i^{k+1}-\y_i^{k+1} \right \rangle+h_i(\z_i^{k+1})\right)+\frac{L_i(\thetak)^2}{2}\|\z_i^{k+1}-\z_i^k\|^2\label{p1l2}\\
&=&(1-\thetak)\left(g_i(\y_i^{k+1})+ \left \langle\nabla g_i(\y_i^{k+1}),\x_i^k-\y_i^{k+1} \right \rangle+h_i(\x_i^k)\right)\notag\\
&&+\thetak\left(g_i(\y_i^{k+1})+ \left \langle\nabla g_i(\y_i^{k+1}),\x_i-\y_i^{k+1} \right \rangle+ \left \langle\nabla g_i(\y_i^{k+1}),\z_i^{k+1}-\x_i \right \rangle+h_i(\z_i^{k+1})\right)\notag\\
&&+\frac{L_i(\thetak)^2}{2}\|\z_i^{k+1}-\z_i^k\|^2\label{p1l3}\\
&\leq&(1-\thetak)f_i(\x_i^k)+\thetak\left(g_i(\x_i)+ \left \langle\nabla g_i(\y_i^{k+1}),\z_i^{k+1}-\x_i \right \rangle+h_i(\z_i^{k+1})\right)+\frac{L_i(\thetak)^2}{2}\|\z_i^{k+1}-\z_i^k\|^2\label{p1l4}\\
&\leq&(1-\thetak)f_i(\x_i^k)+\thetak\left(g_i(\x_i)+h_i(\x_i)+ \left \langle\A_i^T(\bm{\hat\lambda}^{k+1})+(L_i\thetak+\betak\eta_i)(\z_i^{k+1}-\z_i^k),\x_i-\z_i^{k+1} \right \rangle\right)\notag\\
&&+\frac{L_i(\thetak)^2}{2}\|\z_i^{k+1}-\z_i^k\|^2\label{p1l5}\\
&=&(1-\thetak)f_i(\x_i^k)+\thetak f_i(\x_i)+\thetak \left \langle\A_i^T(\bm{\hat\lambda}^{k+1}),\x_i-\z_i^{k+1} \right \rangle\notag\\
&&-\frac{L_i(\thetak)^2+\thetak\betak\eta_i}{2}\left(\|\z_i^{k+1}-\x_i\|^2-\|\z_i^k-\x_i\|^2+\|\z_i^{k+1}-\z_i^k\|^2\right)+\frac{L_i(\thetak)^2}{2}\|\z_i^{k+1}-\z_i^k\|^2\label{p1l6},
\end{eqnarray}
where (\ref{p1l1}) uses (\ref{Lem_lips}), (\ref{pro3566}) uses (\ref{updatexi}), (\ref{p1l2}) is from the convexity of $h_i$, (\ref{p1l4}) is from the convexity of $g_i$, (\ref{p1l5}) uses (\ref{eq:opt_z2}) and (\ref{p1l6}) uses (\ref{lempytha2}). Rerangging the above inequality leads to
\begin{eqnarray}
&&\left(f_i(\x_i^{k+1})-f_i(\x_i)\right)-{\thetak} \left \langle\A_i^T(\hatlambdakk),\x_i-\z_i^{k+1} \right \rangle\notag\\
&\leq&({1-\thetak})\left(f_i(\x_i^k)-f_i(\x_i)\right) + \frac{L_i(\thetak)^2}{2}\left(\|\z_i^{k}-\x_i\|^2-\|\z_i^{k+1}-\x_i\|^2\right)\notag\\
&&+\frac{\thetak\betak\eta_i}{2}\left(\|\z_i^{k}-\x_i\|^2-\|\z_i^{k+1}-\x_i\|^2-\|\z_i^{k+1}-\z_i^k\|^2\right)\label{pro3proof622}
\end{eqnarray}
Dividing both sides of the above inequality by $(\thetak)^2$ leads to
\begin{eqnarray*}
	&&\frac{1}{(\thetak)^2}\left(f_i(\x_i^{k+1})-f_i(\x_i)\right)-\frac{1}{\thetak} \left \langle\A_i^T({\hatlambdakk}),\x_i-\z_i^{k+1} \right \rangle\notag\\
	&\leq&\frac{1-\thetak}{(\thetak)^2}\left(f_i(\x_i^k)-f_i(\x_i)\right)+\frac{L_i}{2}\left(\|\z_i^{k}-\x_i\|^2-\|\z_i^{k+1}-\x_i\|^2\right)\notag\\
	&&+\frac{\betak\eta_i}{2\thetak}\left(\|\z_i^{k}-\x_i\|^2-\|\z_i^{k+1}-\x_i\|^2-\|\z_i^{k+1}-\z_i^k\|^2\right)
\end{eqnarray*}
The proof is completed by using  $\frac{1-\thetakk}{(\thetakk)^2}=\frac{1}{(\thetak)^2}$. $\hfill\blacksquare$

\begin{proposition}\label{prop4}
	In Algorithm \ref{alg_fastplalmps}, the following result holds for any $\blambda$
	\begin{eqnarray}
	&& \langle\A(\zkk)-\b,\blambda-\hatlambdakk\rangle+\frac{\betak\alpha}{2}\|\A(\zkk)-\b\|^2 \notag\\
	& \leq &  \frac{1}{2\betak} \left(\|\blambdak-\blambda\|^2-\|\blambdakk-\blambda\|^2\right)+\frac{\betak}{2}\sumi\eta_i\|\zkk_i-\zk_i\|^2,\label{pro4eq7222333}
	\end{eqnarray}
	where $\alpha=\min\left\{\frac{1}{n+1},\left\{\frac{\eta_i-n\|\A_i\|^2}{(n+1)\|\A_i\|^2},i=1,\cdots,n\right\}\right\}$.
\end{proposition}
\emph{Proof.} By using (\ref{update_lambda}) and (\ref{lempytha2}), we have
\begin{eqnarray}
&& \langle \A(\zkk)-\b,\blambda-\hatlambdakk\rangle \notag\\
&=&  \frac{1}{\betak}\langle\blambda^{k+1}-\blambdak,\blambda-\hatlambdakk\rangle \notag\\
&=&\frac{1}{2\betak}\left(\|\blambda-\blambdak\|^2-\|\blambda-\blambdakk\|^2\right)-\frac{1}{2\betak}\left(\|\hatlambdakk-\blambdak\|^2-\|\blambdakk-\hatlambdakk\|^2\right).\label{pro4eq722}
\end{eqnarray}
Now, consider the last two terms in the above inequality. We deduce
\begin{eqnarray}
&& \frac{1}{2\betak}\left(\|\hatlambdakk-\blambdak\|^2-\|\blambdakk-\hatlambdakk\|^2\right) \notag\\
&=& \frac{\betak}{2}\left(\left\|\sumi\A_i(\zk_i)-\b\right\|^2- \left\|\sumi\A_i(\zkk_i-\zk_i)\right\|^2\right) \notag\\
&\geq&  \frac{\betak}{2}\left(\left\|\sumi\A_i(\zk_i)-\b\right\|^2-\sumi n\|\A_i\|^2\|\zkk_i-\zk_i\|^2 \right) \notag\\
&=& \frac{\betak}{2}\left(\left\|\sumi\A_i(\zk_i)-\b\right\|^2 + \sumi\frac{\eta_i-n\|\A_i\|^2}{\|\A_i\|^2}\|\A_i\|^2\|\zkk_i-\zk_i\|^2-\sumi\eta_i\|\zkk_i-\zk_i\|^2\right) \notag\\
&\geq&\frac{\betak}{2}\left(\alpha (n+1)\left(\left\|\sumi\A_i(\zk_i)-\b\right\|^2 + \sumi\|\A_i(\zkk_i-\zk_i)\|^2 \right)-\sumi\eta_i\|\zkk_i-\zk_i\|^2\right)\notag\\
&\geq&\frac{\betak}{2}\left(\alpha \left\|\sumi\A_i(\zkk_i)-\b\right\|^2-\sumi\eta_i\|\zkk_i-\zk_i\|^2\right)\\
&=& \frac{\betak\alpha}{2}\|\A(\zkk)-\b\|^2 -\frac{\betak}{2}\sumi\eta_i\|\zkk_i-\zk_i\|^2 \label{pro4eq7445}
\end{eqnarray}
The proof is completed by substituting (\ref{pro4eq7445}) into (\ref{pro4eq722}).
$\hfill\blacksquare$

\begin{theorem}\label{con_fastPS}
	In Algorithm \ref{alg_fastplalmps}, for any $K>0$, we have
	\begin{equation}\label{con_FastPS}
	\begin{split}
	&  f(\x^{K+1})-f(\x^*)+\left \langle\bm{\lambda}^*,\A(\x^{K+1})-\b \right \rangle+\frac{\beta\alpha}{2}\left\|\A(\x^{K+1})-\b\right\|^2 \\
	\leq \ & \frac{2L_{\max}D^2_{\x^*}}{(K+2)^2}+\frac{2\beta\eta_{\max} D^2_{X}}{K+2}+\frac{2D^2_{\Lambda}}{\beta(K+2)},
	\end{split}
	\end{equation}
	where $\alpha=\min\left\{\frac{1}{n+1},\left\{\frac{\eta_i-n\|\A_i\|^2}{2(n+1)\|\A_i\|^2},i=1,\cdots,n\right\}\right\}$,  $L_{\max}=\max\{L_i,i=1,\cdots,n\}$ and $\eta_{\max}=\max\{\eta_i,i=1,\cdots,n\}$.
\end{theorem}


\emph{Proof.} Let $\x_i=\x_i^*$ and $\blambda=\blambda^*$ in (\ref{pro331})  and (\ref{pro4eq7222333}). We have
\begin{eqnarray}
&&\frac{1-\thetakk}{(\thetakk)^2}\sum_{i=1}^n\left(f_i(\x_i^{k+1})-f_i(\x_i^*)\right)-\frac{1-\thetak}{(\thetak)^2}\sum_{i=1}^n\left(f_i(\x_i^k)-f_i(\x_i^*)\right)+\frac{1}{\thetak}\sumi \left \langle\bm\lambda^*,\A_i(\z_i^{k+1})-\b \right \rangle \label{thm2prof7888}\\
&\leq&\frac{1}{2}\sum_{i=1}^nL_i\left(\|\z_i^{k}-\x_i^*\|^2-\|\z_i^{k+1}-\x_i^*\|^2\right)+\frac{\betak}{2\thetak}\sum_{i=1}^n\eta_i\left(\|\z_i^{k}-\x_i^*\|^2-\|\z_i^{k+1}-\x_i^*\|^2-\|\zkk_i-\zk_i\|^2\right)\notag\\
&&\frac{1}{\thetak}\langle\blambda^*-\hatlambdakk,\A(\zkk)-\b\rangle \label{thm2prof7999999}\\
&\leq&\frac{1}{2}\sum_{i=1}^nL_i\left(\|\z_i^{k}-\x_i^*\|^2-\|\z_i^{k+1}-\x_i^*\|^2\right)+\frac{\betak}{2\thetak}\sum_{i=1}^n\eta_i\left(\|\z_i^{k}-\x_i^*\|^2-\|\z_i^{k+1}-\x_i^*\|^2-\|\zkk_i-\zk_i\|^2\right)\notag\\
&&+\frac{1}{2\thetak\betak}\left(\|\bm{\lambda}^k-\bm\lambda^*\|^2-\|\bm{\lambda}^{k+1}-\bm\lambda^*\|^2 \right)+\frac{\betak}{2\thetak}\sumi\eta_i\|\zkk_i-\zk_i\|^2-\frac{\betak\alpha}{2\thetak}\|\A(\zkk)-\b\|^2\label{thm2proof80000}\\
&=&\frac{1}{2}\sum_{i=1}^nL_i\left(\|\z_i^{k}-\x_i^*\|^2-\|\z_i^{k+1}-\x_i^*\|^2\right)+\frac{\betak}{2\thetak}\sum_{i=1}^n\eta_i\left(\|\z_i^{k}-\x_i^*\|^2-\|\z_i^{k+1}-\x_i^*\|^2\right)\notag\\
&&+\frac{1}{2\thetak\betak}\left(\|\bm{\lambda}^k-\bm\lambda^*\|^2-\|\bm{\lambda}^{k+1}-\bm\lambda^*\|^2 \right)-\frac{\betak\alpha}{2\thetak}\|\A(\zkk)-\b\|^2\label{pro2proof8111}
\end{eqnarray}
where (\ref{thm2prof7999999}) uses the fact $\A(\x^*)=\b$ and (\ref{thm2proof80000}) uses (\ref{pro4eq7222333}).

Summing (\ref{thm2prof7888})-(\ref{pro2proof8111}) from $k=0$ to $K$ and fixing $\betak=\beta>0$, we have
\begin{eqnarray}
&&\frac{1-\theta^{(K+1)}}{(\theta^{(K+1)})^2}\sum_{i=1}^n\left(f_i(\x_i^{K+1})-f_i(\x_i^*)\right)-\frac{1-\theta^{(0)}}{(\theta^{(0)})^2}\sum_{i=1}^n\left(f_i(\x_i^0)-f_i(\x_i^*)\right)+\sumk\frac{1}{\thetak} \left \langle\bm\lambda^*,\A(\z^{k+1})-\b \right \rangle \label{thm4profeq82} \\
&\leq&\frac{1}{2}\sum_{i=1}^nL_i\|\z_i^{0}-\x_i^*\|^2+\sum_{k=0}^K\frac{1}{2\thetak}\sum_{i=1}^n\beta\eta_i\left(\|\z_i^{k}-\x_i^*\|^2-\|\z_i^{k+1}-\x_i^*\|^2\right)\notag\\
&&+\sum_{k=0}^K\frac{1}{2\thetak\beta}\left (\|\bm{\lambda}^k-\bm\lambda^*\|^2-\|\bm{\lambda}^{k+1}-\bm\lambda^*\|^2 \right)-\sumk\frac{\beta\alpha}{2\thetak}\|\A(\zkk)-\b\|^2\notag\\ &\leq&\frac{1}{2}\sum_{i=1}^nL_i\|\z_i^{0}-\x_i^*\|^2+\frac{1}{2}\sum_{k=0}^K\sum_{i=1}^n\beta\eta_i\|\z_i^{k}-\x_i^*\|^2+\frac{1}{2\beta}\sum_{k=0}^K\|\bm{\lambda}^k-\bm\lambda^*\|^2-\sumk\frac{\beta\alpha}{2\thetak}\|\A(\zkk)-\b\|^2,\label{pro4proof84445}\\
&\leq&\frac{1}{2}\left(L_{\max}D^2_{\x^*}+K\beta\eta_{\max}D^2_{\X}+\frac{1}{\beta}KD_{\bm\Lambda}-\beta\alpha\sumk\frac{1}{\theta^{(k)}}\|\A(\zkk)-\b\|^2\right)\label{thm2proofproof84565}
\end{eqnarray}
where (\ref{pro4proof84445}) uses (\ref{lemseq333}). Also note that $\theta^{(0)}=1$. So the second term of
(\ref{thm4profeq82}) disappears.

Note that (\ref{pro1proofeq46666}) and (\ref{pro1proofeq499}) also holds here. Substituting (\ref{pro1proofeq46666}) into (\ref{thm4profeq82}) and (\ref{pro1proofeq499}) into (\ref{thm2proofproof84565}) respectively and using $\frac{1-\theta^{(K+1)}}{(\theta^{(K+1)})^2}=\frac{1}{(\theta^{(K)})^2}$, we obtain
\begin{eqnarray*}
	&&\frac{1}{(\theta^{(K)})^2}\sum_{i=1}^n\left(f_i(\x_i^{K+1})-f_i(\x_i^*)\right)+\frac{1}{(\theta^{(K)})^2} \left \langle\bm\lambda^*,\A(\x^{K+1})-\b \right \rangle+\frac{\beta\alpha}{2(\theta^{(K)})^2}\|\A(\x^{K+1})-\b\|^2 \\
	&\leq&\frac{1}{2}\left(L_{\max}D^2_{\x^*}+K\beta\eta_{\max}D^2_{X}+\frac{1}{\beta}KD_{\Lambda}\right).
\end{eqnarray*}
The proof is completed by multiplying both sides of the above inequality with $\theta^{(K)}$ and using (\ref{lem3pro3}).
$\hfill\blacksquare$

%
%

\begin{theorem}\label{bound_lambda}
	Assume the mapping $\A(\x_1,\cdots,\x_n)=\sum_{i=1}^n\A_i(\x_i)$ is onto\footnote{This assumption is equivalent to that the matrix $A\equiv(A_1,\cdots,A_n)$ is of full row rank, where $A_i$ is the matrix representation of $\A_i$.}, the sequence $\{\z^k\}$ is bounded, $\partial h(\x)$ and $\nabla g(\x)$ are bounded if $\x$ is bounded, then $\{\x^k\}$, $\{\y^k\}$ and $\{\blambda^k\}$ are bounded.
\end{theorem}
\emph{Proof.}
Assume $\|\z^k\|\leq C_1$ for all $k$ and $\|\x^0\|\leq C_1$. Then from (\ref{updatexi}) we can easily get $\|\x^k\|\leq C_1$ for all $k$. Then from (\ref{updateyykkk}) we have $\|\y^k\|\leq C_1$ for all $k$. Assume $\|\partial h(\x)\|\leq C_2$ and $\|\nabla g(\x)\|\leq C_2$ if $\|\x\|\leq C_1$. Then from (\ref{pro353}), we have
\begin{eqnarray}
0&\in& \partial h(\z^{k+1})+\nabla g(\y^{k+1})+\A^T(\bm{\lambda}^k)+\betak\A^T(\A(\z^{k})-\b)+\left[\begin{matrix}
(L_1\thetak+\betak\eta_1)(\z_1^{k+1}-\z_1^k) \\
\vdots\\
(L_i\thetak+\betak\eta_i)(\z_i^{k+1}-\z_i^k) \\
\vdots\\
(L_n\thetak+\betak\eta_n)(\z_n^{k+1}-\z_n^k)\notag
\end{matrix}\right],
\end{eqnarray}
and
\begin{eqnarray}
-\A\A^T(\bm{\lambda}^k) &\in& \A\left(\partial h(\z^{k+1})+\nabla g(\y^{k+1})+\betak\A^T(\A(\z^{k})-\b)+\left[\begin{matrix}
(L_1\thetak+\betak\eta_1)(\z_1^{k+1}-\z_1^k) \\
\vdots\\
(L_i\thetak+\betak\eta_i)(\z_i^{k+1}-\z_i^k) \\
\vdots\\
(L_n\thetak+\betak\eta_n)(\z_n^{k+1}-\z_n^k)\notag
\end{matrix}\right]\right).
\end{eqnarray}
So we have
\begin{eqnarray}
\|\bm{\lambda}^k\| &\leq& \left\|(\A\A^T)^{-1}\A\left(\partial h(\z^{k+1})+\nabla g(\y^{k+1})+\betak\A^T(\A\z^{k}-\b)+\left[\begin{matrix}
(L_1\thetak+\betak\eta_1)(\z_1^{k+1}-\z_1^k) \notag\\
\vdots\\
(L_i\thetak+\betak\eta_i)(\z_i^{k+1}-\z_i^k) \\
\vdots\\
(L_n\thetak+\betak\eta_n)(\z_n^{k+1}-\z_n^k)
\end{matrix}\right]\right)\right\|\\
&\leq& \|(\A\A^T)^{-1}\A\|\left(\|\partial h(\z^{k+1})\|+\|\nabla g(\y^{k+1})\|+\|\betak\A^T\A\z^{k}\|+\|\betak\A^T\b\|+(L_{\max}\thetak+\betak\eta_{\max})(\|\z^{k+1}\|+\|\z^k\|)\right)\notag\\
&\leq& \|(\A\A^T)^{-1}\A\|\left( 2C_2+\betak\|\A^T\A\|C_1+\betak\|\A^T\b\|+2(L_{\max}\thetak+\betak\eta_{\max})C_1 \right)\notag.
\end{eqnarray}
for all $k$.
$\hfill\blacksquare$

\end{document}